\documentclass{amsart}
\usepackage{amssymb, latexsym}
\usepackage{mathtools,color}
\usepackage{enumerate}
\usepackage{footmisc}
\usepackage{framed}

\newcommand\dive{\operatorname{div}}

\newcommand \lesssimo{\lesssim_{ \| u_{1} \|_{H^{s-1}}, \,  \| u_{0} \|_{H^{s}}}}

\theoremstyle{plain}
\newtheorem{thm}{Theorem}

\newtheorem{rem}[thm]{Remark}
\newtheorem{prop}[thm]{Proposition}
\newtheorem{lem}[thm]{Lemma}

\begin{document}

\title[Radial $3D$ NLKG]{Introduction to scattering for radial $3D$  NLKG below energy norm}
\author{Tristan Roy}
\address{Nagoya math department, Japan}
\email{tristanroy@math.nagoya-u.ac.jp}

\begin{abstract}
We prove scattering for some radial $3D$ semilinear Klein-Gordon
equations with rough data. First we prove Strichartz-type estimates
in mixed norm spaces. Then by using these decays we establish some
local bounds. By combining these results with a Morawetz-type estimate
and a radial Sobolev inequality we control the variation of an
almost conserved quantity on arbitrarily large intervals. Once we have
showed that this quantity is controlled, we prove that some of these
local bounds can be upgraded to global bounds. This is enough to
establish scattering. All the estimates involved require a delicate
analysis due to the nature of the nonlinearity and the lack of
scaling.
\end{abstract}

\maketitle

\vspace{-0.3in}

\section{Introduction}

In this paper we consider the $p$- defocusing Klein-Gordon equation
on $\mathbb{R}^{3}$

\begin{equation}
\begin{array}{ll}
\partial_{tt} u  - \Delta u + u & =  -|u|^{p-1}u \\
\end{array}
\label{Eqn:NlkgWdat}
\end{equation}
with data $u(0)=u_{0}$, $\partial_{t}u(0)=u_{1}$  lying in $H^{s}$,
$H^{s-1}$ respectively. Here $H^{s}$ is the standard inhomogeneous
Sobolev space i.e $H^{s}$ is the completion of the Schwartz space
$\mathcal{S}(\mathbb{R}^{3})$ with respect to the norm

\begin{equation}
\begin{array}{ll}
\| f \|_{H^{s}} & :=  \| \langle D \rangle^{s} f
\|_{L^{2}(\mathbb{R}^{3})}
\end{array}
\end{equation}
where $\langle D \rangle $ is the operator defined by

\begin{equation}
\begin{array}{ll}
\widehat{\langle D \rangle ^{s} f}(\xi) & := ( 1 + |\xi|^{2})^{\frac{s}{2}}
\hat{f}(\xi)
\end{array}
\end{equation}
and $\hat{f}$ denotes the Fourier transform
\begin{equation}
\begin{array}{ll}
\hat{f}(\xi) & := \int_{\mathbb{R}^{3}} f(x) e^{-i x \cdot \xi} \,
dx
\end{array}
\end{equation}
We are interested in the strong solutions of the $p$- defocusing
Klein-Gordon equation on some interval $[0,T]$ i.e maps $u$,
$\partial_{t} u$ that lie in $C \left([0, \, T], \,
H^{s}(\mathbb{R}^{3}) \right)$, $C \left( [0, \,T], \, H^{s-1} (
\mathbb{R}^{3}) \right)$ respectively and that satisfy

\begin{equation}
\begin{array}{ll}
u(t) & = \cos{(t \langle D \rangle )} u_{0} + \frac{\sin(t \langle D
\rangle)}{\langle D \rangle} u_{1} - \int_{0}^{t} \frac{\sin \left(
(t-t^{'}) \langle D \rangle \right)}{\langle D \rangle} \left(
|u|^{p-1}(t^{'}) u(t^{'}) \right) \, dt^{'}
\end{array}
\label{Eqn:StrongSol}
\end{equation}
The $p$- defocusing Klein-Gordon equation is closely related to the
$p$- defocusing wave equation i.e

\begin{equation}
\begin{array}{ll}
\partial_{tt} v - \triangle v = & - |v|^{p-1} v
\end{array}
\label{Eqn:DefocPWaveEq}
\end{equation}
with data $v(0)=v_{0}$, $\partial_{t} v(0)=v_{1}$.
(\ref{Eqn:DefocPWaveEq}) enjoys the following scaling property

\begin{equation}
\begin{array}{ll}
v(t,x) & \rightarrow \frac{1}{\lambda^{\frac{2}{p-1}}} u \left( \frac{t}{\lambda}, \frac{x}{\lambda} \right) \\
v_{0}(x) & \rightarrow \frac{1}{\lambda^{\frac{2}{p-1}}} u_{0}  \left( \frac{x}{\lambda} \right) \\
v_{1}(x) & \rightarrow \frac{1}{\lambda^{\frac{2}{p-1}+1}} u_{1}
\left( \frac{x}{\lambda} \right)
\end{array}
\label{Eqn:WaveScaling}
\end{equation}
We define the critical exponent $s_{c}:= \frac{3}{2} -
\frac{2}{p-1}$. One can check that the $\dot{H}^{s_{c}} \times
\dot{H}^{s_{c}-1}$ norm of $(u_{0},u_{1})$ is invariant under the
transformation (\ref{Eqn:WaveScaling}) \footnote{ Here $\dot{H}^{m}$
denotes the standard homogeneous Sobolev space endowed with the norm
$\| f \|_{\dot{H}^{m}}: = \| D^{m} f \|_{L^{2}(\mathbb{R}^{3})}$}.
(\ref{Eqn:DefocPWaveEq}) was demonstrated to be locally well-posed
by Lindblad and Sogge \cite{linsog} in $H^{s} \times H^{s-1}$, $s
>\frac{3}{2} - \frac{2}{p-1}$, $p>3$ by using an iterative argument.
In fact their results extend immediately to (\ref{Eqn:NlkgWdat})
\footnote{by rewriting for example (\ref{Eqn:NlkgWdat}) in the
"wave" form $\partial_{tt}u - \triangle u= -|u|^{p-1} u -u$}.

If $p=5$ then $s_{c}=1$ and this is why we say that that the
nonlinearity $|u|^{p-1} u$ is $\dot{H}^{1}$ critical. If $3<p<5$
then $s_{c}<1$ and the regime is $\dot{H}^{1}$ subcritical.

It is well-known that smooth solutions to (\ref{Eqn:NlkgWdat}) have
a conserved energy

\begin{equation}
\begin{array}{ll}
E(u(t)) & := \frac{1}{2} \int_{\mathbb{R}^{3}} \left| \partial_{t} u
(t,x) \right|^{2} \, dx  + \frac{1}{2} \int_{\mathbb{R}^{3}} |
\nabla  u
(t,x)|^{2} \, dx + \frac{1}{2} \int_{\mathbb{R}^{3}} | u(t,x) |^{2} \, dx \\
& + \frac{1}{p+1} \int_{\mathbb{R}^{3}} |u(t,x)|^{p+1} \, dx
\end{array}
\label{Eqn:DefEnergy}
\end{equation}
In fact by standard limit arguments the energy conservation law
remains true for solutions $(u,\partial_{t}u) \in H^{s} \times
H^{s-1}$, $s \geq 1$.

Since the lifespan of the local solution depends only on the $H^{s}
\times H^{s-1}$ norm of the initial data $(u_{0},u_{1})$ (see
\cite{linsog}) then it suffices to find an a priori pointwise in
time bound in $H^{s} \times H^{s-1}$ of the solution
$(u,\partial_{t}u)$ to establish global well-posedness. The energy
captures the evolution in time of the $H^{1} \times L^{2}$ norm of
the solution. Since it is conserved we have global existence of
(\ref{Eqn:NlkgWdat}).

The scattering theory (namely, the existence of the bijective wave
operators) in the energy space  \footnote{i.e with data
$(u_{0},u_{1}) \in H^{1} \times L^{2}$} for (\ref{Eqn:NlkgWdat}) has
been extensively studied for a large range of exponents $p$. In
particular Brenner \cite{bren,brenscat} was able to prove that if $
\frac{7}{3} < p < 5 $, then every solution scatters as $T$ goes to
infinity. In fact he showed scattering for all dimension $n$, $n
\geq 3$  and for all exponent $p$ that is $\dot{H}^{1}$ subcritical
and $L^{2}$ supercritical \footnote{since if $p > 1 + \frac{4}{n}$
then $s_{c}>0$}, i.e $ 1 + \frac{4}{n} < p < 1 +\frac{4}{n-2} $.
Later Nakanishi (\cite{nakanash1}, \cite{nakanash2}) was able to
extend these results to $n=1$ and $2$.

In this paper we are interested in proving scattering results for
data below the energy norm i.e for $s<1$. We will assume that
(\ref{Eqn:NlkgWdat}) has radial data. The main result of this paper
is the following one

\begin{thm}
The $p$-radial defocusing Klein-Gordon equation on $\mathbb{R}^{3}$
is globally well-posed in $ H^{s} \times H^{s-1} $, $ 1 > s > s(p)$
and there exists a scattering state $\left( u_{+,0}, u_{+,1} \right)
\in H^{s} \times H^{s-1}$ such that

\begin{equation}
\begin{array}{ll}
\lim \limits_{T \rightarrow \infty} \| \left( u(T) , \partial_{t}
u(T) \right) - K(T) ( u_{+,0}, \,  u_{+,1} ) \|_{H^{s} \times H^{s-1}}
& = 0
\end{array}
\label{Eqn:ScatteringRes}
\end{equation}
with

\begin{equation}
K(t):= \left(
\begin{array}{lll}
\cos{(t \langle D \rangle)} & & \frac{ \sin {(t \langle D
\rangle)}}{\langle D \rangle} \\
- \langle D \rangle \sin { t \langle D \rangle} & &  \cos{ (t
\langle D \rangle)}
\end{array}
\right)
\end{equation}
$3<p<5$ and

\begin{equation}
\begin{array}{ll}
s_{p} & := \left\{
\begin{array}{l}
1 - \frac{(5-p)(p-3)}{2(p-1)(p-2)}, \,  3< p \leq 4 \\
1 - \frac{(5-p)^{2}}{2(p-1)(6-p)}, \, 4 \leq p < 5
\end{array}
\right.
\end{array}
\end{equation}
\label{Thm:Scat}
\end{thm}

Throughout the paper $\nabla$ denotes the gradient operator. Let
$s_{c}$, $\theta_{1}$,...,$\theta_{3}$ denote the following numbers

\begin{equation}
\begin{array}{ll}
s_{c} & : = \frac{3}{2} - \frac{2}{p-1}
\end{array}
\end{equation}

\begin{equation}
\begin{array}{ll}
\theta_{1} & := \left\{
\begin{array}{l}
\frac{(2s-1)(4-p)}{s(p-1)(p-2)}, \,  3< p \leq 4 \\
\frac{(4s-1)(p-4)}{s(p-1)(6-p)}, \, 4 \leq p < 5
\end{array}
\right.
\end{array}
\end{equation}

\begin{equation}
\begin{array}{ll}
\theta_{2} & := \left\{
\begin{array}{l}
\frac{(p+2)(p-3)}{(p-1)(p-2)}, \,   3 < p \leq 4 \\
\frac{(p+2)(5-p)}{(6-p)(p-1)}, \, 4 \leq p < 5
\end{array}
\right.
\end{array}
\label{Eqn:DefTheta2}
\end{equation}
and

\begin{equation}
\begin{array}{ll}
\theta_{3} & := \left\{
\begin{array}{l}
 \frac{4-p}{s(p-1)(p-2)}, \,   3 < p \leq 4 \\
\frac{p-4}{s(p-1)(6-p)}, \,   4 \leq p < 5
\end{array}
\right.
\end{array}
\end{equation}
We write $F(v)$ for the following function

\begin{equation}
\begin{array}{ll}
F(v) & := |v|^{p-1} v
\end{array}
\end{equation}

Let $I$ be the following multiplier

\begin{equation}
\begin{array}{ll}
\widehat{If}(\xi) & := m(\xi) \hat{f}(\xi)
\end{array}
\end{equation}
where $m(\xi): =  \eta \left( \frac{\xi}{N} \right)$, $\eta$ is a
smooth, radial, nonincreasing in $|\xi|$ such that

\begin{equation}
\begin{array}{ll}
\eta (\xi) & :=  \left\{
\begin{array}{l}
1, \, |\xi| \leq 1 \\
\left( \frac{1}{|\xi|} \right)^{1-s}, \, |\xi| \geq 2
\end{array}
\right.
\end{array}
\end{equation}
and $N \gg 1$ is a dyadic number playing the role of a parameter to be
chosen. We shall abuse the notation and write $m (|\xi|)$ for
$m(\xi)$, thus for instance $m(N)=1$.

Some estimates that we establish throughout the paper require a
Paley-Littlewood decomposition. We set it up now. Let $\phi(\xi)$ be
a real, radial, nonincreasing function that is equal to $1$ on the
unit ball $\left\{ \xi \in \mathbb{R}^{3}: \, |\xi| \leq 1 \right\}$
and that that is supported on $\left\{ \xi \in \mathbb{R}^{3}: \,
|\xi| \leq 2 \right\}$. Let $\psi$ denote the function

\begin{equation}
\begin{array}{ll}
\psi(\xi) & := \phi(\xi) - \phi(2 \xi)
\end{array}
\end{equation}
If $(M,M_{1},M_{2}) \in 2^{\mathbb{Z}}$ are dyadic numbers such that
$M_{2}>M_{1}$ we define the Paley-Littlewood operators in the
Fourier domain by

\begin{equation}
\begin{array}{ll}
\widehat{P_{\leq M} f}(\xi) & := \phi \left( \frac{\xi}{M} \right)
\hat{f}(\xi) \\
\widehat{P_{M} f}(\xi) & := \psi \left( \frac{\xi}{M} \right)
\hat{f}(\xi) \\
\widehat{P_{> M} f}(\xi) & := \hat{f}(\xi) - \widehat{P_{\leq M} f}(\xi) \\
\widehat{P_{\ll M } f}(\xi) & : = \widehat{P_{ \leq \frac{M}{128} }f}(\xi) \\
\widehat{P_{ \gtrsim  M } f}(\xi) & : = \widehat{P_{ > \frac{M}{128}}f}(\xi) \\
P_{M_{1} \leq . \leq M_{2}} f  &  := P_{ > M_{2}} f -  P_{ < M_{1}} f \\
\end{array}
\end{equation}
Since $\sum_{M \in  2^{\mathbb{Z}}} \psi \left( \frac{\xi}{M}
\right)=1$ we have

\begin{equation}
\begin{array}{ll}
f & = \sum_{M \in 2^{\mathbb{Z}}} P_{M} f
\end{array}
\end{equation}
Notice also that

\begin{equation}
\begin{array}{ll}
f & = P_{\ll M}f + P_{\gtrsim M} f
\end{array}
\end{equation}
It $T$ is a multiplier with nonnegative symbol $m$  then
$T^{\frac{1}{2}}$ denotes then multiplier with symbol
$m^{\frac{1}{2}}$. For instance $\widehat{P_{M}^{\frac{1}{2}}
f}(\xi)= \psi^{\frac{1}{2}} \left( \frac{\xi}{M} \right)
\widehat{f}(\xi)$.

Throughout this paper we constantly use Strichartz-type estimates .
Notice that some Strichartz estimates for the Klein-Gordon equation
already exist in Besov spaces \cite{ginebvelo}. Here we have chosen
to work in the $L_{t}^{q} L_{x}^{r}$ spaces in order to avoid too
many technicalities. The following proposition is proved in Section
\ref{Sec:ProofStrKg}

\begin{prop}{\textbf{"Strichartz estimates for Klein-Gordon equations in $L_{t}^{q} L_{x}^{r}$ spaces"}}
Assume that $u$ satisfies the following Klein-Gordon equation on
$\mathbb{R}^{d}$, $d \geq 3$

\begin{equation}
\left\{
\begin{array}{ccl}
\partial_{tt} u  - \Delta u + u & = & Q \\
u(0,x)& = & u_{0}(x)    \\
\partial_{t} u(0,x) & = & u_{1}(x)
\end{array}
\right. \label{Eqn:KlGH}
\end{equation}
Let $T \geq 0$. Then

\begin{equation}
\begin{array}{l}
\| u \|_{L_{t}^{q}([0,T]) L_{x}^{r}} + \| \partial_{t} \langle D
\rangle^{-1} u \|_{L_{t}^{q}([0, \,T]) L_{x}^{r}} +  \| u
\|_{L_{t}^{\infty} \left( [0, \, T], \,
H^{m} \right)} +  \| \partial_{t} u \|_{L_{t}^{\infty} \left( [0, \, T], \, H^{m-1}  \right)} \\ \\
\lesssim \| u_{0} \|_{H^{m}} + \| u_{1} \|_{H^{m-1}} + \| Q
\|_{L_{t}^{\tilde{q}}([0,T]) L_{x}^{\tilde{r}}}
\end{array}
\label{Eqn:StrNlkg}
\end{equation}
under the following assumptions

\begin{itemize}

\item $(q,r)$ is $m$- wave admissible, i.e $(q,r)$  lies in the set $\mathcal{W}$ of wave-admissible points

\begin{equation}
\begin{array}{ll}
\mathcal{W} & := \left\{ (q,r): (q,r) \in (2, \infty] \times [2,
\infty), \, \frac{1}{q} + \frac{d-1}{2r} \leq \frac{d-1}{4}
\right\}
\end{array}
\label{Eqn:Waveadm}
\end{equation}
it obeys the following constraint

\begin{equation}
\begin{array}{ll}
\frac{1}{q} + \frac{d}{r} & =\frac{d}{2} -m
\end{array}
\label{Eqn:qrscal}
\end{equation}
and

\begin{equation}
\begin{array}{ll}
(q,r) & \neq \left( 2, \, \frac{2(d-1)}{d-3}  \right)
\end{array}
\label{Eqn:EdPt}
\end{equation}

\item $(\tilde{q},\tilde{r})$ lies in the dual set $\widetilde{\mathcal{W}}$ of $\mathcal{W}$ i.e

\begin{equation}
\begin{array}{ll}
\widetilde{\mathcal{W}} & := \left\{ (\tilde{q},\tilde{r}):
\frac{1}{\tilde{q}} + \frac{1}{q}=1, \, \frac{1}{\tilde{r}} +
\frac{1}{r}=1 \right\}
\end{array}
\end{equation}
and it satisfies the following inequality

\begin{equation}
\begin{array}{ll}
\frac{1}{\tilde{q}} + \frac{d}{\tilde{r}} -2   & = \frac{1}{q} +
\frac{d}{r}
\end{array}
\label{Eqn:ScalingInhomqtildeq}
\end{equation}

\end{itemize}

\label{prop:StrEstLtqLxr}
\end{prop}

\begin{rem}
Notice that the constraints that $(q,r,\tilde{q},\tilde{r})$ must
satisfy are essentially the same to those in the Strichartz
estimates for the wave equation \cite{linsog}. These similarities
are not that surprising. Indeed the relevant operator is  $e^{i t
\langle D \rangle}$, $e^{itD}$ for the Klein-Gordon, wave equations
respectively \footnote{with $D$ multiplier defined by
$\widehat{Df}(\xi):=|\xi| \widehat{f}(\xi)$}. They are similar to
each other on high frequencies.
\end{rem}

Now we explain the main ideas of this paper.

Our first objective is to establish global well-posedness of
(\ref{Eqn:NlkgWdat}) for data in $H^{s} \times H^{s-1}$,$1>s>s(p)$ \footnote{ \label{Gwp} Notice that the
global well-posedness was already studied in \cite{bodaomiao}. Since it is a prerequisite to study scattering, we mention it.}.
Unfortunately since the solution lies in $H^{s} \times H^{s-1}$ pointwise in time and the energy (\ref{Eqn:DefEnergy}) is infinite.
Therefore we introduce the following mollified energy

\begin{equation}
\begin{array}{ll}
E \left( Iu(t) \right) & := \frac{1}{2} \int_{\mathbb{R}^{3}} \left|
\partial_{t} I u(t,x) \right|^{2} \, dx + \frac{1}{2} \int_{\mathbb{R}^{3}} | D I u(t,x)|^{2} \, dx \\
& + \frac{1}{2} \int_{\mathbb{R}^{3}} |Iu(t,x)|^{2} \, dx   +
\frac{1}{p+1} \int_{\mathbb{R}^{3}} |I u(t,x)|^{p+1} \, dx
\end{array}
\label{Eqn:SmthNrj}
\end{equation}
This is the $I$-method originally designed by J. Colliander, M.
Keel, G. Staffilani, H. Takaoka and T. Tao \cite{almckstt} to study
global existence for rough solutions of semilinear Schr\"odinger
equations. Since the multiplier gets closer  to the identity
operator as the parameter $N$ goes to infinity \footnote{formally
speaking} we expect the variation of the smoothed energy  to
approach zero as $N$ grows. However it is not equal to zero and it
needs to be controlled on an arbitrarily large interval. The
semilinear Schr\"odinger and Wave equations have a scaling property.
In \cite{almckstt,visanzhang} the authors were able after scaling to
make the mollified energy at time zero smaller than one. Then by
using the Strichartz estimates they locally bounded some numbers
that allowed them to find an upper bound of its local variation.
Iterating the process they managed to yield an upper bound
\footnote{depending on $N$, the time and the initial data} of its
total variation. Choosing appropriately the parameter $N$ they
bounded it by a constant. Unfortunately the $p$-defocusing
Klein-Gordon equation does not have any scaling symmetry. We need to
control the variation of (\ref{Eqn:SmthNrj}) by a fixed quantity. A
natural choice is a constant $C>1$ multiplied by the mollified
energy $ E(Iu_{0}) :=E(Iu(0))$ at time zero. It occurs that this is
possible if $E(Iu_{0})$ is bounded by a constant depending on the
parameter $N$: see (\ref{Eqn:UnifProc21}) and
(\ref{Eqn:UnifProc22}). But Proposition \ref{Prop:EstInitMolNrj}
shows that $E(Iu_{0})$ is bounded by a power of $N$. Therefore we
can choose $N$ to control the mollified energy as long as $s >s(p)$.
Since the pointwise in time $H^{s} \times H^{s-1}$ norm of the
solution is bounded by the mollified energy (see
(\ref{Eqn:BounduTMolNrj})) we have global well-posedness. \\ Now we
are interested in proving asymptotic completeness by using the
$I$-method. Notice that this method has already been used in
\cite{visanzhang} to prove scattering below the energy norm for
semilinear Schr\"odinger equations with a power type nonlinearity.
We would like to establish (\ref{Eqn:ScatteringRes}). Notice first
that if this result is true then it implies that the pointwise in
time $H^{s} \times H^{s-1}$ bound of the norm of the solution is bounded by a
function that does not depend on time. Therefore in view of the
previous paragraph, the variation of the smoothed energy should not
depend on time $T$. To this end we use some tools. Recall that this
variation is estimated by using local bounds of some quantities,
namely some $Z_{m,s}$ (see Proposition \ref{prop:LocalBd}). We
divide the whole interval $[0,T]$ into subintervals where the
$L_{t}^{p+2} L_{x}^{p+2}$ of $Iu$ is small and we control these
numbers on them by the Strichartz estimates and a continuity
argument. Notice that in this process we are not allowed to create
powers of time $T$ \footnote{by using H\"older locally in time}
since it will eventually force us to choose $N$ as a function of
$T$. We also need to control the $L_{t}^{p+2} L_{x}^{p+2}$ norm of
the solution on $[0,T]$. Morawetz and Strauss \cite{moraw,morstr}
proved a weighted long time estimate ( see (\ref{Eqn:Ms4}))
depending on the energy. Combining this result with a radial Sobolev
inequality (see (\ref{Eqn:RadSobolev})) \footnote{this is the only
place where we rely crucially on the assumption of spherical
symmetry} we can control the $L_{t}^{p+2} L_{x}^{p+2}$ norm of $u$
by some power of the energy. Of course since the solution lies in
$H^{s} \times H^{s-1}$, $s<1$ we cannot use this inequality as such.
Instead we prove an almost Morawetz-Strauss estimate (see
Proposition \ref{Prop:AlmMor} and Proposition \ref{Prop:EstInt} ) by
substituting $u$ for $Iu$ in the establishment of (\ref{Eqn:Ms4}).
This approach was already used in \cite{triroyrad}. Notice here that
the upper bound of (\ref{Eqn:NormIuLt}) does not depend on $T$
either. The almost conservation law (see Proposition \ref{Prop:Acl})
is proved in Section \ref{Sec:AlmConsLawProof} by performing a
low-high frequency decomposition and using the smoothness of $F$
\footnote{ namely $F$  is $C^{1}$ if $p>3$} when we estimate the low
frequency part of the variation. Combining all these tools we are
able to iterate and globally bound the mollified energy and the
$L_{t}^{p+2} L_{x}^{p+2}$ norm of $u$ by a function of $N$ and the
data. These global results allow us to update a local control of the
$Z_{m,s}$ to a global one. It occurs that scattering holds if some
integrals are finite. By using the global control of the $Z_{m,s}$
in the Cauchy criterion we prove these facts. This is enough to
establish scattering.

$\textbf{Acknowledgements}:$ The author would like to thank Terence
Tao for suggesting him  this problem.

\section{Proof of Theorem \ref{Thm:Scat}}

In this section we prove Theorem \ref{Thm:Scat} assuming that the
following propositions are true.

\begin{prop}{"\textbf{Mollified energy at time $0$ is bounded by $N^{2(1-s)}$}"}
Assume that $s_{c} < s < 1$. Then

\begin{equation}
\begin{array}{ll}
E(Iu_{0}) & \lesssim N^{2(1-s)} \left( \| u_{0} \|^{2}_{H^{s}} + \|
u_{1} \|^{2}_{H^{s-1}} + \| u_{0} \|^{p+1}_{H^{s}} \right)
\end{array}
\end{equation}
\label{Prop:EstInitMolNrj}
\end{prop}

\begin{prop}{\textbf{"Local Boundedness"}}
Assume that $u$ satisfies (\ref{Eqn:NlkgWdat}). Let
$\mathcal{M}=[0,s] \cup \{ 1- \} $. There exists $N= N(\| u_{0}
\|_{H^{s}}, \| u_{1} \|_{H^{s-1}}) \gg 1$ such that if $J$, time
interval, satisfies

\begin{equation}
\begin{array}{ll}
\sup_{t \in J} E(Iu(t)) & \leq 3 E(Iu_{0})
\end{array}
\label{Eqn:InducNrj}
\end{equation}
and

\begin{equation}
\begin{array}{ll}
\| I u \|_{L_{t}^{p+2}(J) L_{x}^{p+2}} & \leq   \frac{1}{
N^{+}\left( E(Iu_{0}) \right)^{\frac{1- \theta_{2}}{2 \theta_{2}}}}
\end{array}
\label{Eqn:LocLong}
\end{equation}
then

\begin{equation}
\begin{array}{ll}
Z(J,u) & \lesssim E^{\frac{1}{2}}(Iu_{0})
\end{array}
\label{Eqn:BoundZ}
\end{equation}
where, given a function $v$,

\begin{equation}
\begin{array}{ll}
Z(J,v) & := \sup_{m \in \mathcal{M}} Z_{m,s}(J,v)
\end{array}
\end{equation}
and

\begin{equation}
\begin{array}{ll}
Z_{m,s}(J,v) & := \sup_{ (q,r)-m \, wave \, adm}  \| \partial_{t}
\langle D \rangle ^{-m} I v \|_{L_{t}^{q}(J) L_{x}^{r}} + \| \langle
D \rangle^{1-m} I v \|_{L_{t}^{q}(J) L_{x}^{r}}
\end{array}
\end{equation}
We recall that, if $a$ and $M$ are two real number, then
$M^{a+}:=M^{a+ \alpha}$ and $M^{a-}:=M^{a - \alpha}$ for $ 0 < \alpha \ll
1$. \label{prop:LocalBd}
\end{prop}

\begin{prop}{\textbf{"Almost Conservation Law "}}
Assume that $u$ satisfies (\ref{Eqn:NlkgWdat}). Let $J=[a,b]$ be a
time interval. Let $3 \leq p < 5$ and $s \geq \frac{3p-5}{2p}$.
Then

\begin{equation}
\begin{array}{ll}
\left| \sup_{t \in J} E(Iu(t)) - E(Iu(a)) \right| & \lesssim
\frac{Z^{p+1}(J,u)}{N^{\frac{5-p}{2}-}}
\end{array}
\label{Eqn:Acl}
\end{equation}

 \label{Prop:Acl}
\end{prop}

\begin{rem}
Notice that if $p=3$ then the upper bound is $O \left(
\frac{1}{N^{1-}} \right)$ modulo $Z^{p+1}(J,u)$. This result has
already been established in \cite{triroyrad} for a slighly different
problem, i.e the defocusing cubic wave equation by using a
multilinear analysis.
\end{rem}

\begin{prop}{\textbf{"Estimate of integrals"}}
Let $J$ be a time interval. Let $v$ be a function. Then for $i=1,2$
we have

\begin{equation}
\begin{array}{ll}
| R_{i}(J,v) | & \lesssim \frac{Z^{p+1}(J,v)}{N^{\frac{5-p}{2}-}}
\end{array}
\label{Eqn:BdRi}
\end{equation}
with

\begin{equation}
\begin{array}{ll}
R_{1}(J,v) : = \int_{J} \int_{\mathbb{R}^{3}} \frac{\nabla
Iv(t,x).x}{|x|} \left( F(Iv) - I F(v) \right) \, dx dt
\end{array}
\end{equation}
and

\begin{equation}
\begin{array}{ll}
R_{2}(J,v) := \int_{J} \int_{\mathbb{R}^{3}} \frac{Iv(t,x)}{|x|}
\left( F(Iv) - I F(v) \right) \, dx dt
\end{array}
\end{equation}
\label{Prop:EstInt}
\end{prop}

\begin{prop}{\textbf{"Almost Morawetz-Strauss Estimate"}}
Let $u$ be a solution of (\ref{Eqn:NlkgWdat}) and let $T \geq 0$.
Then
\begin{equation}
\begin{array}{ll}
\int_{0}^{T} \int_{\mathbb{R}^{3}} \frac{|Iu(t,x)|^{p+1}}{|x|} \, dx
dt & \lesssim \sup_{t \in [0,  T]} E(Iu(t)) + R_{1}([0, \, T],u) +
R_{2}([0, \, T],u)
\end{array}
\label{Eqn:AlmMorStrEst}
\end{equation}
\label{Prop:AlmMor}
\end{prop}

These propositions will be proved in the next sections. The proof of
Theorem \ref{Thm:Scat} is made of four steps

\begin{itemize}

\item \textit{Boundedness of the mollified energy and the quantity $ \| I u \|_{L_{t}^{p+2} L_{x}^{p+2}}$}. We will prove that we can control
the mollified energy $E(Iu)$ and the $L_{t}^{p+2} L_{x}^{p+2}$ norm
of $Iu$  on arbitrarily large intervals $[0, T]$, $T \geq 0$. More
precisely let

\begin{equation}
\begin{array}{l}
F_{T} : = \left\{ T^{'} \in [0, \,T]:
\begin{array}{l}
\sup_{t  \in [0, \, T^{'}]} E(Iu(t)) \leq 2E(Iu_{0}), \\
\| I u \|^{p+2}_{L_{t}^{p+2}([0, \, T^{'}]) L_{x}^{p+2}} \leq C
E^{\frac{3}{2}} (I u_{0})
\end{array}
\right\}
\end{array}
\end{equation}
We claim that $F_{T}=[0, \, T]$ for some universal constant $C \geq
0$ and $N=N(\| u_{0} \|_{H^{s}}, \| u_{1} \|_{H^{s-1}}) \gg 1$ to be
chosen later. Indeed

\begin{itemize}

\item $F_{T} \neq \emptyset $ since $0 \in F_{T}$

\item $F_{T}$ is closed by continuity

\item $F_{T}$ is open. Let $\widetilde{T'} \in F_{T}$. By continuity there exists $ \delta > 0$ such that for all $T^{'} \in
(\widetilde{T^{'}} - \delta, \widetilde{T^{'}} + \delta) \cap [0,
T]$ we have

\begin{equation}
\begin{array}{ll}
\sup_{t \in [0, T^{'}]} E(Iu(t)) & \leq 3 E(Iu_{0})
\end{array}
\label{Eqn:InducNrjThm}
\end{equation}
and

\begin{equation}
\begin{array}{ll}
\| I u \|^{p+2}_{L_{t}^{p+2}([0, \, T^{'}]) L_{x}^{p+2}} & \leq 2 C
E^{\frac{3}{2}}(Iu_{0})
\end{array}
\label{Eqn:InducLtp2Lxp2}
\end{equation}
Let $\mathcal{P}=(J_{j})_{1 \leq j \leq l }$ be a partition of $[0,
\,  T^{'}]$ such that $\| I u \|_{L_{t}^{p+2}(J_{j}) L_{x}^{p+2}} =
\frac{1}{N^{+}  E^{\frac{1- \theta_{2}}{2 \theta_{2}}}(Iu_{0})} $
for all $j=1,..,l-1$ and $\| I u \|_{L_{t}^{p+2}(J_{l}) L_{x}^{p+2}}
\leq \frac{1}{N^{+}  E^{\frac{1- \theta_{2}}{2 \theta_{2}}}(Iu_{0})
}$ with $N^{+}$ defined in Proposition \ref{prop:LocalBd}. Then by
(\ref{Eqn:InducLtp2Lxp2})

\begin{equation}
\begin{array}{ll}
l & \lesssim E^{ \frac{(p+2)(1-\theta_{2})}{2 \theta_{2}}
+\frac{3}{2}}  (I u_{0}) N^{+}
\end{array}
\label{Eqn:Estl1}
\end{equation}
By Proposition \ref{prop:LocalBd} and \ref{Prop:Acl} we get after
iteration

\begin{equation}
\begin{array}{ll}
\sup_{t \in [0, \,T]} E(Iu(t)) - E(Iu_{0})  & \lesssim \frac{E^{
\frac{(p+2)(1-\theta_{2})}{2 \theta_{2}} + \frac{3}{2}   +
\frac{p+1}{2}}(Iu_{0}) }{N^{\frac{5-p}{2}-}}
\end{array}
\label{Eqn:MolNrjVar}
\end{equation}
We recall that, if $A$ and $B$ are two real numbers such that $A
\lesssim B$, then the constant determined by $\lesssim$ in $A
\lesssim B$ is the smallest constant among the K such that $A \leq K
B$.

Let $C_{1}$ be the constant determined by $\lesssim$ in
(\ref{Eqn:MolNrjVar}). If we can choose $N \gg 1$ such that

\begin{equation}
\begin{array}{ll}
C_{1} \frac{E^{  \frac{(p+2)(1-\theta_{2})}{2 \theta_{2}}
+\frac{3}{2}  + \frac{p+1}{2}}(Iu_{0}) }   {N^{\frac{5-p}{2}-}} \leq
E(Iu_{0})
\end{array}
\label{Eqn:UnifProc1}
\end{equation}
then $ \sup_{t \in [0, T^{'}]} E(Iu(t)) \leq 2E(Iu_{0})$. The
constraint (\ref{Eqn:UnifProc1}) is equivalent to

\begin{equation}
\begin{array}{ll}
E(I u_{0}) & \leq \frac{ N^{\frac{(5-p)(p-3)}{(p-1)(p-2)}-}}{
C_{1}^{\frac{2(p-3)}{(p-1)(p-2)}}}
\end{array}
\label{Eqn:UnifProc21}
\end{equation}
if $3 < p \leq 4$ and

\begin{equation}
\begin{array}{ll}
E(I u_{0}) & \leq
\frac{N^{\frac{(5-p)^{2}}{(6-p)(p-1)}-}}{C_{1}^{\frac{2(5-p)}{(6-p)(p-1)}}}
\end{array}
\label{Eqn:UnifProc22}
\end{equation}
if $ 4 \leq p < 5$ after plugging (\ref{Eqn:DefTheta2}) into
(\ref{Eqn:UnifProc1}). By Proposition \ref{Prop:EstInitMolNrj} it
suffices to prove that there exists $N= N (\| u_{0} \|_{H^{s}}, \|
u_{1} \|_{H^{s-1}}) \gg 1$ such that

\begin{equation}
\begin{array}{ll}
N^{2(1-s)} \max(\| u_{0} \|_{H^{s}}, \, \| u_{1} \|_{H^{s-1}}, \, \|
u_{0} \|_{H^{s}}^{p+1} ) & \lesssim
N^{\frac{(5-p)(p-3)}{(p-1)(p-2)}-}
\end{array}
\label{Eqn:UnifProc3}
\end{equation}
in order to satisfy (\ref{Eqn:UnifProc21}) and

\begin{equation}
\begin{array}{ll}
N^{2(1-s)} \max(\| u_{0} \|_{H^{s}}, \, \| u_{1} \|_{H^{s-1}}, \, \|
u_{0} \|_{H^{s}}^{p+1} ) & \lesssim
N^{\frac{(5-p)^{2}}{(6-p)(p-1)}-}
\end{array}
\label{Eqn:UnifProc4}
\end{equation}
in order to satisfy (\ref{Eqn:UnifProc22}). Such a choice is
possible if and only if $s> s(p)$. By Proposition \ref{Prop:AlmMor},
Proposition \ref{Prop:EstInt} and (\ref{Eqn:InducNrjThm}) we get
\begin{equation}
\begin{array}{ll}
\int_{0}^{T^{'}} \int_{\mathbb{R}^{3}} \frac{|Iu(t,x)|^{p+1}}{|x|}
\, dx dt & \lesssim E(Iu_{0}) + \frac{E^{
\frac{(p+2)(1-\theta_{2})}{2
\theta_{2}} + \frac{3}{2}   + \frac{p+1}{2}}(Iu_{0}) }{N^{\frac{5-p}{2}-}} \\
& \lesssim E(I u_{0})
\end{array}
\label{Eqn:AlmMorWe}
\end{equation}
Combining (\ref{Eqn:AlmMorWe}) with the well-known pointwise radial
Sobolev inequality

\begin{equation}
\begin{array}{ll}
|Iu(t,x)| & \lesssim  \frac{ \| I u(t,.) \|_{H^{1}}}{|x|}
\end{array}
\label{Eqn:RadSobolev}
\end{equation}
we have

\begin{equation}
\begin{array}{ll}
\| I u \|^{p+2}_{L_{t}^{p+2}([0, T^{'}]) L_{x}^{p+2}} & \lesssim
E^{\frac{3}{2}}(Iu_{0})
\end{array}
\label{Eqn:NormIuLt}
\end{equation}
and we assign to $C$ the constant determined by $\lesssim$ in
(\ref{Eqn:NormIuLt}).

\end{itemize}

\item \emph{Global existence} We have just proved that

\begin{equation}
\begin{array}{ll}
\sup_{t \in [0,T]} E(Iu(t)) & \leq 2 E(Iu_{0})
\end{array}
\label{Eqn:BdNrjNrjIn}
\end{equation}
and

\begin{equation}
\begin{array}{ll}
\| I u \|^{p+2}_{L_{t}^{p+2}([0,T]) _{x}^{p+2}} \leq C
E^{\frac{3}{2}}(Iu_{0})
\end{array}
\end{equation}
for some well-chosen $N=N(\| u_{0} \|_{H^{s}}, \| u_{1}
\|_{H^{s-1}}) \gg 1$ and $1>s>s(p)$. Therefore by Proposition
\ref{Prop:EstInitMolNrj}

\begin{equation}
\begin{array}{ll}
\sup_{t \in [0,T]} E(Iu(t)) & \lesssimo 1
\end{array}
\label{Eqn:Nrj0T}
\end{equation}
and

\begin{equation}
\begin{array}{ll}
\| I u \|^{p+2}_{L_{t}^{p+2}([0,T]) L_{x}^{p+2}} \lesssimo 1
\end{array}
\label{Eqn:NrjLgEst0T}
\end{equation}
Here $A \lesssimo B$ means that there exists a constant \\
$K:=K \left( \| u_{0} \|_{H^{s}}, \| u_{1} \|_{H^{s-1}} \right)$
such that $A \leq K B$. Now by Plancherel and (\ref{Eqn:Nrj0T})

\begin{equation}
\begin{array}{ll}
\| ( u(T),\partial_{t} u(T)) \|_{H^{s} \times H^{s-1}} & \lesssim E(Iu(T)) \\
& \lesssimo 1
\end{array}
\label{Eqn:BounduTMolNrj}
\end{equation}
This proves global well-posedness of (\ref{Eqn:NlkgWdat}) with data
$(u_{0},u_{1}) \in H^{s} \times H^{s-1}$, $1>s>s(p)$ \footnote{See footnote \ref{Gwp}}.
Moreover by continuity we have

\begin{equation}
\begin{array}{ll}
\sup_{t \in \mathbb{R}} E(Iu(t)) & \lesssimo 1
\end{array}
\label{Eqn:BdRNrj}
\end{equation}
and

\begin{equation}
\begin{array}{ll}
\| I u \|^{p+2}_{L_{t}^{p+2}(\mathbb{R}) L_{x}^{p+2}} \lesssimo 1
\end{array}
\label{Eqn:BdLgEst}
\end{equation}

\item \textit{Global estimates}

Let $ \mathcal{P}:= (\widetilde{J}_{j}=[a_{j}, \, b_{j}])_{1 \leq j
\leq \widetilde{l}}$  be a partition of $[0, \infty)$ such that

\begin{equation}
\begin{array}{ll}
\| I u \|_{L_{t}^{p+2}(J_{j}) L_{x}^{p+2}} & =  \frac{1}{ N^{+}
\left( E(Iu_{0}) \right)^{\frac{1- \theta_{2}}{2 \theta_{2}}}}
\end{array}
\label{Eqn:LocLongThm}
\end{equation}
(except maybe the last one), with $N^{+}$ defined in Proposition \ref{prop:LocalBd}. Notice that
from Proposition \ref{Prop:EstInitMolNrj} and (\ref{Eqn:BdLgEst})
the number of intervals $\widetilde{l}$ satisfies

\begin{equation}
\begin{array}{ll}
\tilde{l} & \lesssim E^{\frac{(p+2)(1-\theta_{2})}{2
\theta_{2}}}(Iu_{0}) \\ & \lesssimo 1
\end{array}
\label{Eqn:Esttl}
\end{equation}
Moreover by slightly modifying the steps between
(\ref{Eqn:PlugNrjStr}) and (\ref{Eqn:EstZss}) and by
(\ref{Eqn:BdRNrj}) we have

\begin{equation}
\begin{array}{ll}
Z_{s,s}(J_{j},u) & \lesssim  E^{\frac{1}{2}} (Iu(a_{j})) +  C_{1}
Z_{s,s}^{\theta_{3}(p-1)+1}(J_{j},u) + C_{2} Z_{s,s}^{\theta
(p-1) +1}(J_{j},u) \\
& \lesssim E^{\frac{1}{2}}(Iu_{0})  +  C_{1}
Z_{s,s}^{\theta_{3}(p-1)+1}(J_{j},u) + C_{2} Z_{s,s}^{\theta (p-1)
+1}(J_{j},u)
\end{array}
\label{Eqn:EstZss}
\end{equation}
with $C_{1}$, $C_{2}$, and $\theta$ defined in (\ref{Eqn:DfnC1}),
(\ref{Eqn:DfnC2}) and (\ref{Eqn:Dfntheta}) respectively. Even if it
means increasing the value of $N=N(\| u_{0} \|_{H^{s}}, \| u_{1}
\|_{H^{s-1}}) \gg 1$ in (\ref{Eqn:UnifProc3}) and
(\ref{Eqn:UnifProc4}) we can assume that (\ref{Eqn:Const1}) and
(\ref{Eqn:Const2}) hold. Therefore by Lemma \ref{lem:LgCont} and
Proposition \ref{Prop:EstInitMolNrj} we have

\begin{equation}
\begin{array}{ll}
Z_{s,s}(J_{j},u) & \lesssim E^{\frac{1}{2}}(Iu_{0}) \\
& \lesssimo 1
\end{array}
\label{Eqn:EstZj}
\end{equation}
By (\ref{Eqn:EstZj}) and (\ref{Eqn:Esttl}) we have

\begin{equation}
\begin{array}{ll}
Z_{s,s}(\mathbb{R},u) &  \lesssimo  1
\end{array}
\label{Eqn:EstZjTot}
\end{equation}

\item \textit{Scattering}

Let

\begin{equation}
\begin{array}{ll}
\mathbf{v}(t) & := \left(
\begin{array}{l}
u(t) \\
\partial_{t} u(t)
\end{array}
\right)
\end{array}
\end{equation}

\begin{equation}
\begin{array}{ll}
\mathbf{v}_{0} : = \left(
\begin{array}{l}
u_{0} \\
u_{1}
\end{array}
\right)
\end{array}
\end{equation}
and

\begin{equation}
\begin{array}{ll}
\mathbf{u_{nl}}(t) & = \left(
\begin{array}{l}
\int_{0}^{t} \frac{\sin{\left((t-t^{'}) \langle D \rangle \right)}}{\langle D \rangle} \left( |u|^{p-1}(t^{'}) u(t^{'}) \right) \, d t^{'} \\ \\
-\int_{0}^{t} \cos{ \left( (t-t^{'}) \langle D \rangle \right)}
\left( |u|^{p-1}(t^{'}) u(t^{'}) \right) \, d t^{'}
\end{array}
\right)
\end{array}
\end{equation}
Then we get from (\ref{Eqn:StrongSol})

\begin{equation}
\begin{array}{ll}
\mathbf{v}(t) & = \mathbf{K}(t) \mathbf{v_{0}} - \mathbf{u_{nl}}(t)
\end{array}
\end{equation}
Recall that the solution $u$  scatters in $H^{s} \times H^{s-1}$  if
there exists

\begin{equation}
\begin{array}{ll}
\mathbf{v_{+,0}} & :=  \left(
\begin{array}{l}
u_{+,0} \\
u_{+,1}
\end{array}
\right)
\end{array}
\label{Eqn:DefVplusZer}
\end{equation}
such that

\begin{equation}
\left\|   \mathbf{v}(t)
 -\mathbf{K}(t) \mathbf{v_{+,0}} \right\|_{H^{s} \times H^{s-1}}
\end{equation}
has a limit as $t \rightarrow \infty$ and the limit is equal to $0$.
In other words since $\mathbf{K}(t)$ is bounded on $H^{s} \times
H^{s-1}$ it suffices to prove that the quantity

\begin{equation}
\left\| \mathbf{K}^{-1}(t) \mathbf{v}(t) -  \mathbf{v_{+,0}}
\right\|_{H^{s} \times H^{s-1}}
\end{equation}
has a limit as $t \rightarrow \infty$ and the limit is equal to $0$.
A computation shows that

\begin{equation}
\begin{array}{ll}
\mathbf{K}^{-1}(t) & = \left(
\begin{array}{cc}
\cos{(t \langle D \rangle)} & -\frac{\sin{(t \langle D \rangle)}}{\langle D \rangle} \\ & \\
\langle D \rangle \sin {(t \langle D \rangle)} & \cos{(t \langle D \rangle)}
\end{array}
\right)
\end{array}
\end{equation}
But

\begin{equation}
\begin{array}{ll}
\mathbf{K}^{-1}(t) \mathbf{v}(t) &  =  \mathbf{v_{0}} -
\mathbf{K}^{-1}(t) \mathbf{u_{nl}}(t)
\end{array}
\end{equation}
By Proposition \ref{prop:StrEstLtqLxr} (more precisely by dualizing
$\| e^{it \langle D \rangle} f \|_{L_{t}^{\frac{2}{1-s}}L_{x}^{\frac{2}{s}}} \lesssim \| f \|_{H^{1-s}}$)

\begin{equation}
\begin{array}{l}
\| \mathbf{K}^{-1}(t_{1}) \mathbf{u_{nl}}(t_{1}) -
\mathbf{K}^{-1}(t_{2}) \mathbf{u_{nl}}(t_{2}) \|_{H^{s} \times
H^{s-1}} \\
\lesssim \| |u|^{p-1} u \|_{L_{t}^{\frac{2}{1+s}}([t_{1}, \, t_{2}]) L_{x}^{\frac{2}{2-s}}} \\
 \lesssim \| \langle D \rangle^{1-s} I \left(  |u|^{p-1} u  \right)
\|_{L_{t}^{\frac{2}{1+s}}([t_{1}, \, t_{2}]) L_{x}^{\frac{2}{2-s}} }
\end{array}
\label{Eqn:CauchyScat}
\end{equation}
If we let $J:=[t_{1},t_{2}]$ in (\ref{Eqn:PlugNrjStr}) and  follow
the same steps up to (\ref{Eqn:EstZss}) we get from
(\ref{Eqn:BdNrjNrjIn})

\begin{equation}
\begin{array}{ll}
 \| \langle D \rangle^{1-s} I \left(  |u|^{p-1} u  \right) \|_{L_{t}^{\frac{2}{1+s}}([t_{1}, \, t_{2}]) L_{x}^{\frac{2}{2-s}} } & \lesssim
C_{1} Z_{s,s}^{\theta_{3}(p-1)+1}([t_{1}, \, t_{2}],u) \\
& + C_{2} Z_{s,s}^{\theta(p-1)+1}([t_{1}, \, t_{2}],u)
\end{array}
\label{Eqn:Donemins}
\end{equation}
By (\ref{Eqn:EstZjTot}), (\ref{Eqn:CauchyScat}) and
(\ref{Eqn:Donemins})

\begin{equation}
\begin{array}{ll}
\lim \limits_{t_{1} \rightarrow \infty} \| \mathbf{K}^{-1}(t_{1})
\mathbf{u_{nl}}(t_{1}) - \mathbf{K}^{-1}(t_{2})
\mathbf{u_{nl}}(t_{2}) \|_{H^{s} \times H^{s-1}} & = 0
\end{array}
\end{equation}
uniformly in $t_{2}$. This proves that $\mathbf{K}^{-1}(t) v(t)$ has
a limit in $H^{s} \times H^{s-1}$ as $t$ goes to infinity. Moreover

\begin{equation}
\begin{array}{ll}
\lim \limits_{t \rightarrow \infty} \left\|  \mathbf{v}(t) -
\mathbf{K}(t)\mathbf{v_{+,0}}   \right\|_{H^{s}\times H^{s-1}} & = 0
\end{array}
\end{equation}
with $\mathbf{v_{+,0}}$ defined in (\ref{Eqn:DefVplusZer}),

\begin{equation}
\begin{array}{ll}
u_{+,0} & := u_{0} + \int_{0}^{\infty} \frac{\sin{(t^{'} \langle D
\rangle)}}{\langle D \rangle} \left( |u|^{p-1}(t^{'}) u(t^{'})
\right) \, d t^{'}
\end{array}
\end{equation}
and

\begin{equation}
\begin{array}{ll}
u_{+,1} & := u_{1} - \int_{0}^{\infty} \cos{(t^{'} \langle D
\rangle)} \left( |u|^{p-1}(t^{'}) u(t^{'}) \right) \, d t^{'}
\end{array}
\end{equation}

\end{itemize}

\section{Proof of "Mollified energy at time 0 is bounded by $N^{2(1-s)}$"}
\label{Sec:AlmConsLawProof}

In this section we aim at proving Proposition
\ref{Prop:EstInitMolNrj}. By Plancherel we have

\begin{equation}
\begin{array}{ll}
\| I u_{1} \|^{2}_{L^{2}} & \lesssim \int_{|\xi| \leq 2N}
|\widehat{u_{1}}(\xi)|^{2} \, d \xi +
 \int_{|\xi| \geq 2N} \frac{N^{2(1-s)}}{|\xi|^{2(1-s)}}  | \widehat{u_{1}}(\xi) |^{2} \, d\xi  \\
 & \lesssim N^{2(1-s)} \| u_{1} \|_{H^{s-1}}^{2}
\end{array}
\end{equation}
Similarly

\begin{equation}
\begin{array}{ll}
\| \nabla I u_{0}  \|^{2}_{L^{2}} & \lesssim \int_{|\xi| \leq 2N}
|\xi|^{2} |\widehat{u_{0}}(t, \xi)|^{2} \, d \xi +
\int_{|\xi| \geq 2N} |\xi|^{2} \frac{N^{2(1-s)}}{|\xi|^{2(1-s)}} |\widehat{u_{0}}(\xi)|^{2} \, d\xi \\
& \lesssim N^{2(1-s)} \| u_{0} \|_{H^{s}}^{2}
\end{array}
\end{equation}
Moreover by the assumption $s>s_{c}$

\begin{equation}
\begin{array}{ll}
\| I u_{0} \|^{p+1}_{L^{p+1}} & \lesssim \| P_{ \ll N} u_{0} \|^{p+1}_{L^{p+1}} + \| P_{\gtrsim N} I u_{0} \|^{p+1}_{L^{p+1}} \\
& \lesssim N^{ (p+1) \left( \frac{3(p-1)}{2(p+1)} - s \right)} \| u_{0} \|^{p+1}_{H^{s}} \\
& \lesssim N^{2(1-s)} \| u_{0} \|_{H^{s}}^{p+1}
\end{array}
\end{equation}

\section{Proof of "Local Boundedness"}

Before attacking the proof of Proposition \ref{prop:LocalBd} let us
prove a short lemma

\begin{lem}
Let $x(t)$ be a nonnegative continuous function of time $t$ such
that $x(0)=0$. Let $X$ be a positive constant and let $\alpha_{i},
\, C_{i}$, $i \in \{ 1,..,m \}$ be nonnegative constants such that

\begin{equation}
\begin{array}{ll}
C_{i} X^{\alpha_{i}-1} & \ll 1
\end{array}
\end{equation}
and

\begin{equation}
\begin{array}{ll}
x(t) & \lesssim X +  \displaystyle {\sum_{i=1}^{m}} C_{i}
x^{\alpha_{i}}(t)
\end{array}
\end{equation}
Then

\begin{equation}
\begin{array}{ll}
x(t) & \lesssim X
\end{array}
\label{Eqn:ResLem}
\end{equation}
\label{lem:LgCont}
\end{lem}

\newcommand{\bproof}{\noindent{\bf Proof: }}
If we let $\overline{x}(t):=\frac{x(t)}{X}$ then we have

\begin{equation}
\begin{array}{ll}
\overline{x}(t) & \lesssim 1 + \displaystyle{ \sum_{i=1}^{m}} C_{i}
X^{\alpha_{i}-1} \overline{x}^{\alpha_{i}}(t)
\end{array}
\end{equation}
and $\overline{x}(0)=0$. Applying a continuity argument to
$\overline{x}$ we have $\overline{x}(t) \lesssim 1$. This implies
(\ref{Eqn:ResLem}).

\newcommand{\eproof}{\hfill $\Box$\\}

Plugging $\langle D \rangle ^{1-m}I$ into (\ref{Eqn:StrNlkg}) we have

\begin{equation}
\begin{array}{ll}
Z_{m,s}(J,u) & \lesssim  E^{\frac{1}{2}}(Iu_{0})  + \| \langle D \rangle^{1-m} I (
|u|^{p-1} u ) \|_{L_{t}^{\frac{2}{1+m}}(J) L_{x}^{\frac{2}{2-m}}}
\end{array}
\label{Eqn:PlugNrjStr}
\end{equation}

There are three cases

\begin{itemize}

\item $m=s$. By (\ref{Eqn:PlugNrjStr}), the fractional Leibnitz rule and H\"older inequality

\begin{equation}
\begin{array}{ll}
Z_{s,s}(J,u) & \lesssim E^{\frac{1}{2}}(Iu_{0}) +  \| \langle D \rangle^{1-s} I u
\|_{L_{t}^{\frac{2}{s}}(J) L_{x}^{\frac{2}{1-s}}}
\| |u|^{p-1} \|_{L_{t}^{2}(J) L_{x}^{2}} \\
& \lesssim E^{\frac{1}{2}}(Iu_{0})  + Z_{s,s}(J,u) \| u \|^{p-1}_{L_{t}^{2(p-1)}(J) L_{x}^{2(p-1)}} \\
& \lesssim   E^{\frac{1}{2}}(Iu_{0}) + Z_{s,s}(J,u) \\
&  \left(  \| P_{ \ll N} u \|^{p-1}_{L_{t}^{2(p-1)}(J) L_{x}^{2(p-1)}}
+ \| P_{\gtrsim N} u
\|^{p-1}_{L_{t}^{2(p-1)}(J) L_{x}^{2(p-1)}} \right) \\
\end{array}
\end{equation}
We are interested in estimating $ \| P_{ \ll N} u
\|^{p-1}_{L_{t}^{2(p-1)} L_{x}^{2(p-1)}} $. There are two cases

\begin{itemize}

\item $3 < p \leq 4$. By interpolation and (\ref{Eqn:LocLong}) we have

\begin{equation}
\begin{array}{l}
\| P_{\ll N} u \|^{p-1}_{L_{t}^{2(p-1)}(J) L_{x}^{2(p-1)}}  \lesssim
\| P_{\ll N} u \|^{\theta_{1}(p-1)}_{L_{t}^{\infty}(J) L_{x}^{2}} \|
P_{\ll N} u \|^{\theta_{2}(p-1)}_{L_{t}^{p+2}(J) L_{x}^{p+2}} \\
\| P_{\ll N} u \|^{\theta_{3} (p-1)}_{L_{t}^{\frac{2}{s}}(J)
L_{x}^{\frac{2}{1-s}}}
 \\
 \lesssim  \| I u \|^{\theta_{1}(p-1)}_{L_{t}^{\infty}(J)
L_{x}^{2}} \| I u \|^{\theta_{2}(p-1)}_{L_{t}^{p+2}(J) L_{x}^{p+2}}
\| \langle D \rangle ^{1-s}
I u \|^{\theta_{3} (p-1)}_{L_{t}^{\frac{2}{s}}(J) L_{x}^{\frac{2}{1-s}}}  \\
\lesssim \frac{E^{\frac{(\theta_{1} + \theta_{2} -1)(p-1)}{2}}(Iu_{0})}{N^{+}} \,  Z_{s,s}^{\theta_{3}(p-1)}(J,u) \\
\lesssim \frac{E^{\frac{(-\theta_{3})(p-1)}{2}}(Iu_{0})}{N^{+}} \,
Z_{s,s}^{\theta_{3}(p-1)}(J,u)
\end{array}
\label{Eqn:LowFreqLocBd3}
\end{equation}

\item $p>4$. By interpolation, Sobolev inequality and (\ref{Eqn:LocLong}) we have

\begin{equation}
\begin{array}{l}
\| P_{\ll N} u \|^{p-1}_{L_{t}^{2(p-1)}(J) L_{x}^{2(p-1)}} \\
\lesssim \left( \| P_{\ll N} u
\|^{\theta_{1}(p-1)}_{L_{t}^{\infty}(J) L_{x}^{6}} \| P_{\ll N} u
\|^{\theta_{2}(p-1)}_{L_{t}^{p+2}(J) L_{x}^{p+2}} \| P_{\ll N} u
\|^{\theta_{3} (p-1)}_{L_{t}^{\frac{2}{s}}(J) L_{x}^{\frac{6}{1-s}}}
\right) \\
\lesssim \left( \| \nabla I u
\|^{\theta_{1}(p-1)}_{L_{t}^{\infty}(J) L_{x}^{2}} \| I u
\|^{\theta_{2}(p-1)}_{L_{t}^{p+2}(J) L_{x}^{p+2}} \| \langle D
\rangle^{1-s}
I u \|^{\theta_{3} (p-1)}_{L_{t}^{\frac{2}{s}}(J) L_{x}^{\frac{2}{1-s}}} \right) \\
\lesssim \frac{E^{\frac{(\theta_{1} + \theta_{2} -1)(p-1)}{2}}(Iu_{0})}{N^{+}} \,  Z_{s,s}^{\theta_{3}(p-1)}(J,u) \\
\lesssim \frac{E^{\frac{(-\theta_{3})(p-1)}{2}}(Iu_{0})}{N^{+}} \,
Z_{s,s}^{\theta_{3}(p-1)}(J,u)
\end{array}
\label{Eqn:LowFreqLocBd5}
\end{equation}
\end{itemize}
Now we estimate  $ \| P_{\gtrsim N} u \|^{p-1}_{L_{t}^{2(p-1)}(J)
L_{x}^{2(p-1)}} $. Let

\begin{equation}
\begin{array}{ll}
\theta:= & \frac{1}{s(p-1)}
\end{array}
\label{Eqn:Dfntheta}
\end{equation}
By interpolation we have

\begin{equation}
\begin{array}{l}
\| P_{\gtrsim N} u \|^{p-1}_{L_{t}^{2(p-1)}(J) L_{x}^{2(p-1)}}
\lesssim \| P_{\gtrsim N} u \|^{\theta(p-1)}_{L_{t}^{\frac{2}{s}}(J)
L_{x}^{\frac{2}{1-s}}} \| P_{\gtrsim N} u \|^{(1-\theta)(p-1)}_{L_{t}^{\infty}(J) L_{x}^{\frac{2 \left( s(p-1) -1 \right)}{2s-1}}} \\
\lesssim \frac{\| \langle D \rangle^{1-s} I u
\|^{\theta(p-1)}_{L_{t}^{\frac{2}{s}}(J)
L_{x}^{\frac{2}{1-s}}}}{N^{(1-s) \theta (p-1)}} \frac{ \| \langle D \rangle I u
\|^{(1- \theta)(p-1)}_{L_{t}^{\infty}(J) L_{x}^{2}}}{N^{(p-1)(1-
\theta)(1-s)} N^{+}} \\ \\
\lesssim E^{ \frac{(1-
\theta)(p-1)}{2}}(Iu_{0}) \frac{Z_{s,s}^{\theta
(p-1)}(J,u)}{N^{(1-s)(p-1)} N^{+}}
\end{array}
\label{Eqn:HighFreqLocBd}
\end{equation}
since $s>s_{c} \geq \frac{1}{p-1}$. Therefore we get from
(\ref{Eqn:PlugNrjStr}), (\ref{Eqn:LowFreqLocBd3}),
(\ref{Eqn:LowFreqLocBd5}) and (\ref{Eqn:HighFreqLocBd})

\begin{equation}
\begin{array}{l}
Z_{m,s}(J,u) \lesssim   E^{\frac{1}{2}} (Iu_{0}) +  C_{1}
Z_{s,s}^{\theta_{3}(p-1)+1}(J,u) + C_{2} Z_{s,s}^{\theta (p-1)
+1}(J,u)
\end{array}
\label{Eqn:EstZss}
\end{equation}
with

\begin{equation}
\begin{array}{ll}
C_{1} & :=  \frac{ E^{\frac{- \theta_{3} (p-1)}{2}}(Iu_{0})} {N^{+}}
\end{array}
\label{Eqn:DfnC1}
\end{equation}
and

\begin{equation}
\begin{array}{ll}
C_{2} & := \frac{E^{ \frac{(1- \theta)(p-1)}{2}}(Iu_{0})
}{N^{(1-s)(p-1)} N^{+}}
\end{array}
\label{Eqn:DfnC2}
\end{equation}
Notice that by Proposition \ref{Prop:EstInitMolNrj}

\begin{equation}
\begin{array}{ll}
C_{1}  E^{\frac{\theta_{3}(p-1)}{2}}( Iu_{0}) & \lesssim \frac{1}{N^{+}} \\
& \ll 1
\end{array}
\label{Eqn:Const1}
\end{equation}
and

\begin{equation}
\begin{array}{ll}
C_{2}  E^{\frac{\theta(p-1)}{2}}(Iu_{0})  & \lesssim \frac{1}{N^{+}} \\
& \ll 1
\end{array}
\label{Eqn:Const2}
\end{equation}
if we choose $N= N(\| u_{0}\|_{H^{s}}, \| u_{1} \|_{H^{s-1}}) \gg 1$.
Applying Lemma \ref{lem:LgCont}, we get

\begin{equation}
\begin{array}{ll}
Z_{s,s}(J,u) & \lesssim  E^{\frac{1}{2}}(Iu_{0})
\end{array}
\label{Eqn:Zmsinf}
\end{equation}

\item $m<s$ Notice that by (\ref{Eqn:LowFreqLocBd3}), (\ref{Eqn:LowFreqLocBd5}), (\ref{Eqn:HighFreqLocBd}), (\ref{Eqn:Const1}), (\ref{Eqn:Const2}) and
(\ref{Eqn:Zmsinf})

\begin{equation}
\begin{array}{ll}
\| u \|^{p-1}_{L_{t}^{2(p-1)}(J) L_{x}^{2(p-1)}} & \lesssim \frac{1}{N^{+}} \\
& \ll 1
\end{array}
\label{Eqn:BdLp}
\end{equation}
Moreover

\begin{equation}
\begin{array}{ll}
Z_{m,s}(J,u) & \lesssim E^{\frac{1}{2}}(Iu_{0}) + \| \langle D \rangle^{1-m} I ( |u|^{p-1} u ) \|_{L_{t}^{\frac{2}{1+m}(J)} L_{x}^{\frac{2}{2-m}}} \\
& \lesssim  E^{\frac{1}{2}}(Iu_{0})  + \| \langle D \rangle ^{1-m} I
u \|_{L_{t}^{\frac{2}{m}}(J) L_{x}^{\frac{2}{1-m}}}
\| u \|^{p-1}_{L_{t}^{2(p-1)}(J) L_{x}^{2(p-1)}} \\
& \lesssim E^{\frac{1}{2}}(Iu_{0})   + Z_{m,s}(J,u) \| u
\|^{p-1}_{L_{t}^{2(p-1)}(J) L_{x}^{2(p-1)}}
\end{array}
\label{Eqn:Bdminfs}
\end{equation}
By (\ref{Eqn:BdLp}) and (\ref{Eqn:Bdminfs}) and Lemma
\ref{lem:LgCont}, we get (\ref{Eqn:BoundZ}).

\item $m=1-=1- \alpha$ with $\alpha$ small. We have

\begin{equation}
\begin{array}{l}
Z_{m,s}(J,u) \lesssim E^{\frac{1}{2}}(Iu_{0})  + \| \langle D \rangle ^{1-(1-)} I \left( |u|^{p-1} u \right) \|_{L_{t}^{1+}(J) L_{x}^{2-}} \\
 \lesssim E^{\frac{1}{2}}(Iu_{0})  + N^{+} \| |u|^{p-1}  u \|_{L_{t}^{1^{+}}(J) L_{x}^{2-}} \\
 \lesssim E^{\frac{1}{2}}(Iu_{0})  + \, N^{+} \| \, |P_{\ll N} u|^{p-1} P_{\ll N} u \|_{L_{t}^{1+}(J)L_{x}^{2-}} + \\
 N^{+} \| \, |P_{\ll N} u|^{p-1} P_{\gtrsim N} u  \|_{L_{t}^{1+}(J)
L_{x}^{2-}}   + \, N^{+} \| \, | P_{\gtrsim N} u |^{p-1} P_{\ll N} u
\|_{L_{t}^{1+}(J) L_{x}^{2-}} \\
+ \, N^{+} \| \,  | P_{\gtrsim N}
u |^{p-1} P_{\gtrsim N} u \|_{L_{t}^{1+}(J) L_{x}^{2-}}
\end{array}
\label{Eqn:1min1}
\end{equation}
But by (\ref{Eqn:BdLp}) we have

\begin{equation}
\begin{array}{l}
 N^{+} \| |P_{\ll N} u|^{p-1} P_{\ll N} u \|_{L_{t}^{1+}(J) L_{x}^{2-}} \lesssim N^{+} \| I u \|_{L_{t}^{2+}(J) L_{x}^{\infty-}} \| I u
\|^{p-1}_{L_{t}^{2(p-1)}(J) L_{x}^{2(p-1)}} \\
\lesssim N^{+} \| \langle D \rangle^{1-(1-)} I u \|_{{L_{t}^{2+}(J) L_{x}^{\infty-}}} \| u \|^{p-1}_{L_{t}^{2(p-1)}(J) L_{x}^{2(p-1)}} \\
\lesssim N^{+} \| u \|^{p-1}_{L_{t}^{2(p-1)}(J) L_{x}^{2(p-1)}}  Z_{1-,s}(J,u) \\
\lesssim \frac{Z_{1-,s}(J,u)}{N^{+}}
\end{array}
\end{equation}
Similarly

\begin{equation}
\begin{array}{ll}
\| N^{+} |P_{ \gtrsim N} u|^{p-1} P_{\ll N} u \|_{L_{t}^{1+}(J)
L_{x}^{2-}} & \lesssim  \frac{Z_{1-,s}(J,u)}{N^{+}}
\end{array}
\label{Eqn:1min2}
\end{equation}
Moreover since $s>\frac{p-3}{2}$

\begin{equation}
\begin{array}{l}
N^{+} \| \, | P_{\ll N N} u |^{p-1} P_{\gtrsim N} u \|_{L_{t}^{1+}
L_{x}^{2-}} \\
\lesssim N^{+} \|  P_{\ll N} u
\|^{p-1}_{L_{t}^{(p-1)+}(J) L_{x}^{\frac{6(p-1)}{p-3}-}} \|
P_{\gtrsim N} u \|_{L_{t}^{\infty}(J) L_{x}^{\frac{6}{6-p}-}} \\
\lesssim N^{+} Z^{p-1}_{1-,s}(J,u) \frac{ \| \langle D \rangle I u
\|_{L_{t}^{\infty}(J) L_{x}^{2}}}{N^{\frac{5-p}{2}-}} \\
\lesssim  \frac{ E^{\frac{1}{2}}(I u_{0})}{N^{\frac{5-p}{2}-}}
Z^{p-1}_{1-,s}(J,u)
\end{array}
\label{Eqn:1min3}
\end{equation}
By Proposition \ref{Prop:EstInitMolNrj} we have for $N=N \left( \|
u_{0} \|_{H^{s}}, \| u_{1} \|_{H^{s-1}} \right) \gg 1 $

\begin{equation}
\begin{array}{l}
N^{+} \| \, | P_{\gtrsim N} u |^{p-1} P_{\gtrsim N} u
\|_{L_{t}^{1+}(J) L_{x}^{2 -}} \lesssim N^{+} \frac{ \| \langle D
\rangle^{1- \left( \frac{2}{p} - \right)} I u \|^{p}_{L_{t}^{p+}(J)
L_{x}^{\frac{2p}{p-2}-}}}{ N^{\frac{5-p}{2}}} \\
\lesssim  \frac{ Z^{p}_{\frac{2}{p}-,s}(J,u) }{ N^{\frac{5-p}{2}-}}
\\ \\
\lesssim \frac{  E^{\frac{p}{2}}(I u_{0}) }{ N^{\frac{5-p}{2}-}} \\
\lesssim E^{\frac{1}{2}}(I u_{0})
\end{array}
\label{Eqn:1min4}
\end{equation}
since $s > \frac{3p-5}{2p}$. Now by (\ref{Eqn:1min1}),
(\ref{Eqn:1min2}), (\ref{Eqn:1min3}) and (\ref{Eqn:1min4})

\begin{equation}
\begin{array}{ll}
Z_{1-,s}(J,u) & \lesssim E^{\frac{1}{2}} (Iu_{0}) +
\frac{Z_{1-,s}(J,u)}{N^{+}} +
 \frac{ E^{\frac{1}{2}} (I u_{0}) }   {N^{\frac{5-p}{2}-}} Z^{p-1}_{1-,s}(J,u)
\end{array}
\label{Eqn:Z1minsEst}
\end{equation}
Let $C_{3} := \frac{  E^{\frac{1}{2}}(I u_{0})
}{N^{\frac{5-p}{2}-}}$. Then by Proposition \ref{Prop:EstInitMolNrj}

\begin{equation}
\begin{array}{ll}
C_{3}  E^{\frac{p-2}{2}}(Iu_{0}) & \ll 1
\end{array}
\label{Eqn:C3}
\end{equation}
and

\begin{equation}
\begin{array}{ll}
\frac{1}{N^{+}} & \ll 1
\end{array}
\label{Eqn:BoundoneN}
\end{equation}
if $N = N \left( \| u_{0} \|_{H^{s}}, \| u_{1} \|_{H^{s-1}} \right)
\gg 1$. From Lemma \ref{lem:LgCont}, (\ref{Eqn:BdLp}),
(\ref{Eqn:Z1minsEst}) (\ref{Eqn:C3}) and (\ref{Eqn:BoundoneN}) we
get (\ref{Eqn:BoundZ}).

\end{itemize}

\vspace{2mm}

\section{Proof of "Almost Morawetz-Strauss estimate"}

In this section we prove Proposition \ref{Prop:AlmMor}.

First we recall the proof of the Morawetz-Strauss estimate  based
upon the important equality \cite{moraw, morstr, strauss}

\begin{equation}
\begin{array}{l}
\Re \left(  \left( \frac{\nabla \bar{u} \cdot x}{|x|} +
\frac{\bar{u}}{|x|}   \right) \left( \partial_{tt} u - \triangle u +
u + |u|^{p-1} u
\right) \right)  = \partial_{t}  \left( \Re \left( \left( \frac{\nabla \bar{u} \cdot x}{|x|} + \frac{\bar{u}}{|x|}  \right) \partial_{t} u \right) \right) \\
\\ + \dive \left(
- \frac{|\partial_{t} u|^{2} x }{2 |x|} - \frac{|u|^{2} 
x }{2|x|^{3}} - \Re \left( \left( \frac{\nabla \bar{u} \cdot x}{|x|}
+ \frac{\overline{u}}{|x|} \right) \nabla u \right)  + \frac{|\nabla
u|^{2} x}{2|x|} + \frac{|u|^{p+1} x}{(p+1)|x|} + \frac{|u|^{2} x}{2}
\right) \\
+ \frac{p-1}{p+1} \frac{|u|^{p+1}}{|x|} + \frac{1}{|x|} \left(
|\nabla u|^{2} -\frac{|\nabla u \cdot x|^{2}}{|x|^{2}}  \right)
\end{array}
\label{Eqn:DiveEq}
\end{equation}
Integrating (\ref{Eqn:DiveEq}) with respect to space and time we
have

\begin{equation}
\begin{array}{l}
\frac{p-1}{p+1} \int_{0}^{T} \int_{\mathbb{R}^{3}} \frac{|u|^{p+1}(t,x)}{|x|} \, dx dt + 2 \pi \int_{0}^{T} |u|^{2}(0,t) \, dt \\
= - \int_{\mathbb{R}^{3}} \Re \left( \left( \frac{\nabla \overline{u}(T,x)
\cdot x}{|x|} + \frac{\overline{u}(T,x)}{|x|}  \right)
\partial_{t} u(T,x) \right) \, dx \\
 + \int_{\mathbb{R}^{3}} \Re  \left( \left( \frac{\nabla \overline{u}(0,x) \cdot x}{|x|} + \frac{\overline{u}(0,x)}{|x|} \right)
\partial_{t} u(0,x) \right) \, dx
\end{array}
\label{Eqn:Ms1}
\end{equation}
if $u$ satisfies (\ref{Eqn:NlkgWdat}). By Cauchy-Schwartz

\begin{equation}
\begin{array}{ll}
\left| \int_{\mathbb{R}^{3}} \Re \left( \left( \frac{\nabla
\overline{u} \cdot x}{|x|} + \frac{\overline{u}}{|x|} \right)
\partial_{t} u(T,x) \right) \, dx \right| & \lesssim
E^{\frac{1}{2}}(u) \left( \int_{\mathbb{R}^{3}} \left| \frac{\nabla
\overline{u}(T,x)\cdot x}{|x|} + \frac{\overline{u}(T,x)}{|x|}
\right|^{2} \, dx \right)^{\frac{1}{2}}
\end{array}
\label{Eqn:CauchSchw}
\end{equation}
After expansion we have

\begin{equation}
\begin{array}{ll}
\int_{\mathbb{R}^{3}} \left| \frac{\nabla \overline{u}(T,x). \cdot
x}{|x|} + \frac{\overline{u}(T,x)}{|x|} \right|^{2} \, dx & =
\int_{\mathbb{R}^{3}} \left| \frac{\nabla u(T,x) \cdot x}{|x|}
\right|^{2} \, dx + 2 \int_{\mathbb{R}^{3}} \frac{\nabla \left(
\frac{|u|^{2}(T,x)}{2} \right) \cdot x}{|x|^{2}} \, dx \\
& + \int_{\mathbb{R}^{3}}  \frac{|u|^{2} (T,x)}{|x|^{2}}  \, dx \\
& = \int_{\mathbb{R}^{3}} \left| \frac{\nabla u(T,x) \cdot x}{|x|} \right|^{2} \, dx  \\
& \lesssim E(u)
\end{array}
\label{Eqn:Ms2}
\end{equation}
Here we used the identity

\begin{equation}
\begin{array}{ll}
\dive \left( \frac{|u|^{2}(T,x) x}  {2 |x|^{2}} \right) & =
\frac{\nabla \left( \frac{ |u|^{2}(T,x)}{2} \right) \cdot
x}{|x|^{2}} + \frac{|u|^{2}(T,x)}{2|x|^{2}}
\end{array}
\end{equation}
Hence we get

\begin{equation}
\begin{array}{ll}
\int_{\mathbb{R}^{3}} \left| \frac{\nabla \overline{u}(T,x) \cdot
x}{|x|} + \frac{\overline{u}(T,x)}{|x|} \right|^{2} \, dx & \lesssim
E(u)
\end{array}
\label{Eqn:Ms2}
\end{equation}
Similarly

\begin{equation}
\begin{array}{ll}
\int_{\mathbb{R}^{3}} \left| \frac{\nabla \overline{u}(0,x) \cdot
x}{|x|} + \frac{\overline{u}(0,x)}{|x|} \right|^{2} \, dx & \lesssim
E(u)
\end{array}
\label{Eqn:Ms3}
\end{equation}
We get from (\ref{Eqn:Ms1}), (\ref{Eqn:CauchSchw}), (\ref{Eqn:Ms2}) and (\ref{Eqn:Ms3}) the
Morawetz-Strauss estimate

\begin{equation}
\begin{array}{ll}
\int_{0}^{T} \int_{\mathbb{R}^{3}} \frac{|u|^{p+1}(t,x)}{|x|} \, dx
dt & \lesssim E(u)
\end{array}
\label{Eqn:Ms4}
\end{equation}

Now we plug the multiplier $I$ into (\ref{Eqn:DiveEq}) and we redo
the computations. We get (\ref{Eqn:AlmMorStrEst}).

\section{Proof of "Almost conservation law" and "Estimate of integrals"}

The proof of Proposition \ref{Prop:Acl}, \ref{Prop:EstInt} relies on
the following lemma

\begin{lem}
Let $G$ such that $\| G \|_{L_{t}^{\infty}(J) L_{x}^{2}} \lesssim
Z(J,v)$. If $s \geq \frac{3p-5}{2p}
> s_{c}$ and $3
\leq p < 5$ then

\begin{equation}
\begin{array}{ll}
 \int_{J} \int_{\mathbb{R}^{3}} \left| G \,  ( F(Iv) - IF(v)) \right|     \, dx  dt  & \lesssim \frac{
 Z^{p+1}(J,v)}{N^{\frac{5-p}{2}-}}
\end{array}
\end{equation}
\label{lem:LemVar}
\end{lem}

\bproof We have

\begin{equation}
\begin{array}{l}
\int_{J} \int_{\mathbb{R}^{3}} \left| G  \, ( F(Iv) - IF(v) ) \right| \, dx dt \\ \\
\lesssim \| G \|_{L_{t}^{\infty}(J) L_{x}^{2}} \| F(Iv) - F(v)
\|_{L_{t}^{1}(J) L_{x}^{2}} + \| G \|_{L_{t}^{\infty}(J) L _{x}^{2}}
\| F(v) - I
F(v) \|_{L_{t}^{1}(J) L_{x}^{2}} \\ \\
\lesssim Z(J,v)  \left( \| F(Iv) - F(v) \|_{L_{t}^{1}(J) L_{x}^{2}}
+ \| F(v) - I F(v) \|_{L_{t}^{1}(J) L_{x}^{2}} \right)
\end{array}
\end{equation}
Let
\begin{equation}
\begin{array}{ll}
X_{1}:= & \| F(Iv) - F(v) \|_{L_{t}^{1}(J) L_{x}^{2}}
\end{array}
\end{equation}
and
\begin{equation}
\begin{array}{ll}
X_{2} := & \| F(v)  - I F(v) \|_{L_{t}^{1}(J) L_{x}^{2}} \\
\end{array}
\end{equation}
We are interested in estimating $X_{1}$. By the fundamental theorem
of calculus we have the pointwise bound

\begin{equation}
\begin{array}{ll}
| F(Iv) -F(v)| & \lesssim  \max{ \left( |Iv|,|v| \right)}^{p-1} |Iv -v| \\
\end{array}
\end{equation}
Plugging this bound into $X_{1}$ we get

\begin{equation}
\begin{array}{l}
X_{1}  \lesssim \| P_{\ll N} v
\|^{p-1}_{L_{t}^{\frac{4(p-1)}{7-p}+}(J)
L_{x}^{\frac{4(p-1)}{p-3}-}} \| P_{\gtrsim N} v
\|_{L_{t}^{\frac{4}{p-3}-}(J) L_{x}^{\frac{4}{5-p}+}} \\
+ \| P_{\gtrsim N} v \|^{p-1}_{L_{t}^{p}(J) L_{x}^{2p}} \|
P_{\gtrsim N} v\|_{L_{t}^{p}(J) L_{x}^{2p}} \\  \lesssim
\frac{1}{N^{\frac{5-p}{2}-}} \left(
\begin{array}{l}
\| \langle D \rangle^{1-(1-)} I v
\|^{p-1}_{L_{t}^{\frac{4(p-1)}{7-p}+}(J)
L_{x}^{\frac{4(p-1)}{p-3}-}} \|
\langle D \rangle^{1-(\frac{p-3}{2}+)} I v  \|_{L_{t}^{\frac{4}{p-3} -}(J) L_{x}^{\frac{4}{5-p}+} }  + \\
\| \langle D \rangle^{1- \frac{3p-5}{2p}} I v \|^{p}_{L_{t}^{p}(J)
L_{x}^{2p}}
\end{array}
\right) \\
\lesssim \frac{Z^{p-1}_{1-,s}(J,v) Z_{\frac{p-3}{2}+,s}(J,v) \, + \,
Z^{p}_{\frac{3p-5}{2p},s}(J,v)}{N^{\frac{5-p}{2}-}} \\
\lesssim \frac{Z^{p}(J,v)}{N^{\frac{5-p}{2}-}}
\end{array}
\end{equation}
Now we turn to  $X_{2}$. On low frequencies we use the smoothness of
$F$ whereas on high frequencies we take advantage of the regularity
of $v$. More precisely by the fundamental theorem
of calculus we have
\begin{equation}
\begin{array}{l}
F(v)  = F( P_{\ll N} v + P_{\gtrsim N} v) \\
= F(P_{\ll N} v) +
\left( \int_{0}^{1} |P_{\ll N} v + s P_{\gtrsim N} v|^{p-1} \, ds
\right) P_{\gtrsim N } v  \\
+ \left( \int_{0}^{1} \frac{P_{\ll N} v + s P_{\gtrsim N} v}{
\vspace{0.2mm} \overline{P_{\ll N} v + s P_{\gtrsim N} v}  } | P_{\ll
N} v + s P_{\gtrsim N} v|^{p-1} \, ds \right) \overline{P_{\gtrsim
N} v}
\end{array}
\end{equation}
Therefore
\begin{equation}
\begin{array}{ll}
X_{2} \lesssim \| P_{\gtrsim N} F(v) \|_{L_{t}^{1}(J) L_{x}^{2}} \\
\lesssim \| P_{\gtrsim N} F(P_{\ll N} v) \|_{L_{t}^{1}(J) L_{x}^{2}}
+ \| \, |P_{\ll N} v|^{p-1} \, P_{\gtrsim N} v \|_{L_{t}^{1}(J)
L_{x}^{2}} \\
+ \| \,  | P_{\gtrsim N} v|^{p-1} \, P_{\gtrsim N} v \|_{L_{t}^{1}(J) L_{x}^{2}} \\
\lesssim X_{2,1} + X_{2,2} + X_{2,3}
\end{array}
\end{equation}
with  $X_{2,1}:=\| P_{\gtrsim N} F(P_{\ll N} v) \|_{L_{t}^{1}(J)
L_{x}^{2}} $, $ X_{2,2}:= \| \, |P_{\ll  N} v|^{p-1} \, P_{\gtrsim N}
v \|_{L_{t}^{1}(J) L_{x}^{2}} $ and $X_{2,3}:= \| \,  | P_{\gtrsim
N} v|^{p-1} \, P_{\gtrsim N} v \|_{L_{t}^{1}(J) L_{x}^{2}} $. But
again by the fundamental theorem of calculus

\begin{equation}
\begin{array}{l}
X_{2,1} \lesssim \frac{1}{N}  \|  \nabla F(P_{\ll N} v)
\|_{L_{t}^{1}(J) L_{x}^{2}} \\
\lesssim \frac{1}{N} \| \, |P_{\ll N} v|^{p-1} \nabla P_{\ll N} v + \frac{ |P_{\ll N} v|^{p-1} P_{\ll N}
v}{ \overline{P_{\ll N} v}} \overline{\nabla P_{\ll N} v}
\|_{L_{t}^{1}(J) L_{x}^{2}}
\end{array}
\end{equation}
Therefore

\begin{equation}
\begin{array}{l}
X_{2,1} \lesssim \frac{1}{N} \| P_{\ll N } v
\|_{L_{t}^{\frac{4(p-1)}{7-p}+}(J)
L_{x}^{\frac{4(p-1)}{p-3}-}}^{p-1}
 \| \nabla P_{\ll N} v \|_{L_{t}^{\frac{4}{p-3}-}(J) L_{x}^{\frac{4}{5-p}+}} \\
\lesssim  \frac{1}{N^{\frac{5-p}{2}-}} \| \langle D \rangle^{1-(1-)} I v
\|_{L_{t}^{\frac{4(p-1)}{7-p}+}(J)
L_{x}^{\frac{4(p-1)}{p-3}-}}^{p-1}
 \| \langle D \rangle^{1- \left( \frac{p-3}{2} + \right)} I v \|_{L_{t}^{\frac{4}{p-3}-}(J) L_{x}^{\frac{4}{5-p}+}} \\
\lesssim \frac{Z^{p-1}_{1-,s}(J,v)
Z_{\frac{p-3}{2}+,s}(J,v)}{N^{\frac{5-p}{2}-}} \\
\lesssim \frac{Z^{p}(J,v)}{N^{\frac{5-p}{2}-}}
\end{array}
\end{equation}
Moreover

\begin{equation}
\begin{array}{l}
X_{2,2} \lesssim \frac{1}{N^{\frac{5-p}{2}-}} \| \langle D \rangle
^{1-(1-)} I v \|_{L_{t}^{\frac{4(p-1)}{7-p}+}(J)
L_{x}^{\frac{4(p-1)}{p-3}-}}^{p-1}  \| \langle D \rangle^{1- \left(
\frac{p-3}{2} + \right)} I v
\|_{L_{t}^{\frac{4}{p-3} -}(J) L_{x}^{\frac{4}{5-p}+}} \\
\lesssim \frac{Z^{p-1}_{1-,s}(J,v)
Z_{\frac{p-3}{2}+,s}(J,v)}{N^{\frac{5-p}{2}-}} \\
\lesssim \frac{Z^{p}(J,v)}{N^{\frac{5-p}{2}-}}
\end{array}
\end{equation}
As for $X_{2,3}$ we have

\begin{equation}
\begin{array}{l}
X_{2,3} \lesssim \| P_{\gtrsim N} v \|^{p}_{L_{t}^{p}(J) L_{x}^{2p}}
\\  \lesssim \frac{ \| \langle D \rangle^{1- \frac{3p-5}{2p} }
I v \|^{p}_{L_{t}^{p}(J) L_{x}^{2p}} } {N^{\frac{5-p}{2}}}  \\
\lesssim \frac{ Z_{\frac{3p-5}{2p},s}^{p} (J,v)} {N^{\frac{5-p}{2}}}
\\ \lesssim \frac{Z^{p}(J,v)}{N^{\frac{5-p}{2}-}}
\end{array}
\end{equation}
\eproof Let $t^{'} \in J=[a,b]$. Then if $u$ is a solution to
(\ref{Eqn:NlkgWdat}) then

\begin{equation}
\begin{array}{l}
\left| E(Iu(t^{'})) - E(Iu(a))  \right| = \left| \int_{[a,t^{'}]}
\int_{\mathbb{R}^{3}} \Re ( \overline{ \partial_{t} I u}
(F(Iu)-IF(u)) ) \right| \\
\lesssim \int_{[a,t^{'}]} \int_{\mathbb{R}^{3}} \left| \overline{
\partial_{t} I u} (F(Iu)-IF(u)) \right|
\end{array}
\label{Eqn:VarNrj}
\end{equation}
Notice that

\begin{equation}
\begin{array}{l}
\| \partial_{t} I u \|_{L_{t}^{\infty}(J) L_{x}^{2}}   \lesssim
Z_{0,s}(J,u)
\end{array}
\end{equation}
Applying Lemma \ref{lem:LemVar}  with $ G:= \partial_{t} I u$ to
(\ref{Eqn:VarNrj}) we get (\ref{Eqn:Acl}). Notice also that

\begin{equation}
\begin{array}{l}
\| \frac{ \nabla I v \cdot x}{|x|} \|_{L_{t}^{\infty}(J) L_{x}^{2}}  \lesssim \| \nabla I v \|_{L_{t}^{\infty}(J) L_{x}^{2}} \\
\lesssim Z_{0,s}(J,v)
\end{array}
\label{Eqn:BdGrad}
\end{equation}
and that

\begin{equation}
\begin{array}{l}
\| \frac{I v}{|x|}  \|_{L_{t}^{\infty}(J) L_{x}^{2}}  \lesssim \| \nabla I v \|_{L_{t}^{\infty}(J) L_{x}^{2}} \\
\lesssim Z_{0,s}(J,v)
\end{array}
\end{equation}
by Hardy inequality. Letting $G(t,x):=\frac{ \nabla I v(t,x) .
x}{|x|}$ we get (\ref{Eqn:BdRi}) from (\ref{Eqn:BdGrad}) and Lemma
\ref{lem:LemVar} for $i=1$. Similarly (\ref{Eqn:BdRi}) holds for
$i=2$ if we let $G(t,x):= \frac{I v(t,x)}{|x|}$.

\vspace{5mm}


\section{Strichartz estimates for $NLKG$ in $L_{t}^{q} L_{x}^{r}$ spaces}
\label{Sec:ProofStrKg}

The techniques used in the proof of these estimates are, broadly
speaking, standard  \cite{linsog,keeltao}. However some subtleties
appear because unlike the homogeneous Schrodinger and wave equations
the homogeneous defocusing Klein-Gordon equation does not enjoy any
scaling property. Now we mention them. Regarding the estimates
involving the homogeneous part of the solution we apply, broadly
speaking, a $" TT^{*} "$ argument to the truncated cone operators
localized at all the frequencies \footnote{i.e to $e^{i t \langle D \rangle}
P_{M}$, $M \in 2^{\mathbb{Z}}$ : see (\ref{Eqn:TMEst})} instead of
applying it at frequency equal to one and then use a scaling
argument for the other frequencies. The inhomogeneous estimates are
slightly more complicated to establish. In the first place we try to
reduce the estimates (see (\ref{Eqn:HghFreqInhom})) localized at all
frequencies to the estimate at frequency one (see
(\ref{Eqn:IneqMInterm})). This strategy does not totally work
because of the lack of scaling. However the remaining estimate (see
\ref{Eqn:RemainInt}), after duality is equivalent to an homogeneous
estimate on high frequencies (see (\ref{Eqn:ChVar})) that has
already been established.

Let $u$ be the solution of (\ref{Eqn:KlGH}) with data $(u_{0},
u_{1})$. We can substitute $[0, \, T]$ for $\mathbb{R}$ in
(\ref{Eqn:StrNlkg}) without loss of generality.

Let  $u_{l}(t):= \cos{(t \langle D \rangle)} u_{0} + \frac{\sin{(t
\langle D \rangle)} }{\langle D \rangle} u_{1}$ and $u_{nl}(t):=
-\int_{0}^{t} \frac{\sin{(t-t^{'}) \langle D \rangle}}{\langle D
\rangle} Q(t^{'})$. We need to show

\begin{equation}
\begin{array}{ll}
\| u_{l} \|_{ L_{t}^{\infty} H^{m}} + \| \partial_{t} u_{l}(t)
\|_{L_{t}^{\infty} H^{m-1}} & \lesssim \| u_{0} \|_{H^{m}} + \|
u_{1} \|_{H^{m-1}}
\end{array}
\label{Eqn:HmLin1}
\end{equation}

\begin{equation}
\begin{array}{ll}
\| u_{l} \|_{L_{t}^{q} L_{x}^{r}} + \| \partial_{t} \langle D
\rangle^{-1} u_{l} \|_{L_{t}^{q} L_{x}^{r}}  & \lesssim \| u_{0}
\|_{H^{m}} + \| u_{1} \|_{H^{m-1}}
\end{array}
\label{Eqn:HmLin2}
\end{equation}

\begin{equation}
\begin{array}{ll}
\| u_{nl} \|_{L_{t}^{q} L_{x}^{r}} + \| \partial_{t} \langle D
\rangle^{-1} u_{nl} \|_{L_{t}^{q}  L_{x}^{r}} & \lesssim \| Q
\|_{L_{t}^{\tilde{q}} L_{x}^{\tilde{r}}}
\end{array}
\label{Eqn:HmNLin1}
\end{equation}
and

\begin{equation}
\begin{array}{ll}
\| u_{nl} \|_{L_{t}^{\infty}  H^{m}} + \| \partial_{t} u_{nl}
\|_{L_{t}^{\infty}  H^{m-1}} & \lesssim \| Q \|_{L_{t}^{\tilde{q}}
L_{x}^{\tilde{r}}}
\end{array}
\label{Eqn:HmNLin2}
\end{equation}
By Plancherel theorem we have (\ref{Eqn:HmLin1}). We prove
(\ref{Eqn:HmLin2}), (\ref{Eqn:HmNLin1}) and (\ref{Eqn:HmNLin2}) in
the next subsections.

\vspace{2mm}

\subsection{Proof of (\ref{Eqn:HmLin2})}

By decomposition and substitution it suffices to prove

\begin{equation}
\begin{array}{ll}
\| e^{it  \langle D \rangle} u_{0} \|_{L_{t}^{q} L_{x}^{r}} &
\lesssim \| u_{0} \|_{H^{m}}
\end{array}
\label{Eqn:HmLin2a}
\end{equation}
If we could prove for every Schwartz function $f$

\begin{equation}
\begin{array}{ll}
\| e^{it  \langle D \rangle} P_{\leq 1} f \|_{L_{t}^{q} L_{x}^{r}} &
\lesssim \| f \|_{L^{2}}
\end{array}
\label{Eqn:StrLocLin1}
\end{equation}
and

\begin{equation}
\begin{array}{ll}
\| e^{i t \langle D \rangle} P_{M} f  \|_{L_{t}^{q} L_{x}^{r}} &  \lesssim M^{m} \| f  \|_{L^{2}} \\
\end{array}
\label{Eqn:StrLocLinM}
\end{equation}
for $M \in 2^{\mathbb{Z}}$, $M > 1$, then (\ref{Eqn:HmLin2}) would
follow. Indeed let $\widetilde{P_{M}}:= P_{ \frac{M}{2} \leq  \,
\leq 2M}$ and $\widetilde{P_{ \leq 1}} := P_{\leq 2}$. Applying
(\ref{Eqn:StrLocLinM}) to $f:= \widetilde{P_{M}} f$ we have

\begin{equation}
\begin{array}{ll}
\| e^{i t \langle D \rangle } P_{M} f  \| _{L_{t}^{q} L_{x}^{r}} & \lesssim M^{m} \| \widetilde{P_{M}} f \|_{L^{2}} \\
& \lesssim \| \widetilde{P_{M}} f  \|_{ \dot{H^{m}}}
\end{array}
\end{equation}
Similarly plugging $f:= \widetilde{P_{ \leq 1}} f$ into
(\ref{Eqn:StrLocLin1}) we have

\begin{equation}
\begin{array}{ll}
\| e^{i t \langle D \rangle } P_{\leq 1} f  \| _{L_{t}^{q} L_{x}^{r}} & \lesssim  \| \widetilde{P_{\leq 1}} f \|_{L^{2}} \\
& \lesssim \| f \|_{H^{m}}
\end{array}
\label{Eqn:LfPaley}
\end{equation}
Before moving forward, we recall the fundamental Paley-Littlewood
equality \cite{stein}: if $1< p <  \infty$ and $h$ is Schwartz then

\begin{equation}
\begin{array}{ll}
\| h \|_{L^{p}} & \sim \| \left( \sum_{M \in 2^{\mathbb{Z}}}  |P_{M}
h|^{2} \right)^{\frac{1}{2}} \|_{L^{p}}
\end{array}
\label{Eqn:PLIneq}
\end{equation}
We plug $h:=P_{>1} f$ into (\ref{Eqn:PLIneq}). Hence by Minkowski
inequality and Plancherel theorem

\begin{equation}
\begin{array}{ll}
\| e^{i t \langle D \rangle} P_{>1} f \|_{L_{t}^{q} L_{x}^{r}} &
\lesssim
\| \left( \sum_{ M \geq 1} | e^{ i t \langle D \rangle } P_{M} f |^{2} \right)^{\frac{1}{2}} \|_{L_{t}^{q} L_{x}^{r} } \\
& \lesssim  \left(  \sum_{M \geq 1} \| e^{ i t \langle D \rangle} P_{M} f \|^{2}_{L_{t}^{q} L_{x}^{r}}  \right)^{\frac{1}{2}} \\
& \lesssim \left( \sum_{M \geq 1} \| \widetilde{P_{M}} f \|^{2}_{\dot{H}^{m}} \right)^{\frac{1}{2}} \\
& \lesssim  \| f \|_{H^{m}}
\end{array}
\label{Eqn:HfPaley}
\end{equation}
since $q \geq 2$ and $r \geq 2$. Combining (\ref{Eqn:LfPaley}) with
(\ref{Eqn:HfPaley}) we get (\ref{Eqn:HmLin2a}). It remains to prove
(\ref{Eqn:StrLocLin1}) and (\ref{Eqn:StrLocLinM}). Let $T_{1}(f):=
e^{it \langle D \rangle} P_{\leq 1} f$ and $T_{M}(f):= e^{ i t
\langle D \rangle} P_{M} f$, $M \in 2^{\mathbb{Z}}$, $M > 1$. We
have

\begin{equation}
\begin{array}{ll}
T_{1}(f)(t,x) & := \int  \phi(\xi) e^{i t \langle
\xi \rangle} \hat{f}(\xi) e^{i \xi.x} \, d \xi
\end{array}
\end{equation}
and if $M \in 2^{\mathbb{Z}}$, $M>1$ let

\begin{equation}
\begin{array}{ll}
T_{M}(f)(t,x) & := \int  \psi \left( \frac{\xi}{M}
\right) e^{i t \langle \xi \rangle } \hat{f}(\xi)  e^{i \xi.x} \, d
\xi
\end{array}
\end{equation}
We would like to prove

\begin{equation}
\begin{array}{ll}
\| T_{1} (f) \|_{L_{t}^{q} L_{x}^{r}} & \lesssim \| f \|_{L^{2}}
\end{array}
\end{equation}
and

\begin{equation}
\begin{array}{ll}
\| T_{M} (f) \|_{L_{t}^{q} L_{x}^{r}} & \lesssim \| f \|_{H^{m}}
\end{array}
\label{Eqn:TMEst}
\end{equation}
By a $"TT^{*}"$ argument we are reduced showing for every continuous
in time Schwartz in space function $g$

\begin{equation}
\begin{array}{ll}
\| T_{1} T_{1}^{*} (g)  \|_{L_{t}^{q} L_{x}^{r}} & \lesssim \| g
\|_{L_{t}^{q^{'}} L_{x}^{r^{'}}}
\end{array}
\label{Eqn:ResT1}
\end{equation}
and similarly

\begin{equation}
\begin{array}{ll}
\| T_{M} T^{*}_{M} (g) \|_{L_{t}^{q}  L_{x}^{r}} & \lesssim M^{2m}
\| g \|_{L_{t}^{q^{'}}  L_{x}^{r^{'}}}
\end{array}
\label{Eqn:ResTM}
\end{equation}
with $\frac{1}{q}+ \frac{1}{q^{'}}=1$ and $\frac{1}{r} +
\frac{1}{r^{'}}=1$. But a computation shows that

\begin{equation}
\begin{array}{ll}
T_{1} T_{1}^{*} (g)  & = K_{1} \ast g \\
& = \int K_{1}(t-t^{'},.) \ast g(t^{'},.) \, dt^{'}
\end{array}
\label{Eqn:kern1}
\end{equation}
and

\begin{equation}
\begin{array}{ll}
T_{M}T_{M}^{*} (g) & = K_{M} \ast g \\
& = \int K_{M}(t-t^{'},.) \ast g(t^{'},.) \, dt^{'}
\end{array}
\label{Eqn:kernM}
\end{equation}
with

\begin{equation}
\begin{array}{ll}
K_{1}( t-t^{'},x ) & : = \int  |\phi (\xi)|^{2} e^{i
\langle \xi \rangle (t-t^{'})} e^{i \xi \cdot x} \, d \xi
\end{array}
\end{equation}
and

\begin{equation}
\begin{array}{ll}
K_{M}(t-t^{'},x) & : = \int \left| \psi \left(
\frac{\xi}{M}  \right) \right|^{2} e^{i  \langle \xi \rangle
(t-t^{'})} e^{i \xi \cdot x } \, d \xi
\end{array}
\end{equation}
One one hand by Plancherel equality we have

\begin{equation}
\begin{array}{ll}
\| K_{M}(t-t^{'},.) \ast g(t^{'},.) \|_{L^{2}} & \lesssim \|
g(t^{'},.) \|_{L^{2}}
\end{array}
\label{Eqn:PlK1}
\end{equation}
On the other hand

\begin{equation}
\begin{array}{ll}
\| K_{M}(t-t^{'},.) \ast g(t^{'},.)  \|_{L^{\infty}} & \lesssim  \|
K_{M}(t-t^{'},.) \|_{L^{\infty}} \| g(t^{'},.) \|_{L^{1}}
\end{array}
\end{equation}
where $\| K_{M}(t-t^{'},.) \|_{L^{\infty}}$ is estimated by the
stationary phase method  \cite{ginebvelo}, p 441

\begin{equation}
\begin{array}{ll}
\| K_{M} (t-t^{'},.) \|_{L^{\infty}} & \lesssim M^{d} \min{ \left(1,
\frac{1}{\left( M |t-t^{'}| \right)^{\frac{d-1}{2}}} \right)} \min{
\left( 1, \left( \frac{M}{|t-t^{'}|} \right)^{\frac{1}{2}} \right)}
\end{array}
\end{equation}
and

\begin{equation}
\begin{array}{ll}
\| K_{1} (t-t^{'},.)  \|_{L^{\infty}} & \lesssim \min { \left(1,
\frac{1}{|t-t^{'}|^{\frac{d}{2}}} \right) }
\end{array}
\end{equation}
By complex interpolation we have

\begin{equation}
\begin{array}{ll}
\| K_{1}(t-t^{'},.) \ast g(t^{'},.) \|_{L^{r}} & \lesssim  \left(
\min { \left(1, \frac{1}{|t-t^{'}| ^{\frac{d}{2}}} \right) }
\right)^{1 - \frac{2}{r}} \| g (t^{'},.)\|_{L^{r{'}}}
\end{array}
\label{Eqn:EstK1}
\end{equation}
and

\begin{equation}
\begin{array}{ll}
\| K_{M}(t-t^{'},.) \ast g(t^{'},.) \|_{L^{r}} & \lesssim
\widetilde{K_{M}}(t-t^{'})  \| g(t^{'},.) \|_{L^{r^{'}}}
\end{array}
\label{Eqn:IneqConv}
\end{equation}
with

\begin{equation}
\begin{array}{ll}
\widetilde{K_{M}}(t)  & := \left( M^{d} \min{ \left(1,
\frac{1}{\left( M |t| \right)^{\frac{d-1}{2}}} \right)} \min{ \left(
1, \left( \frac{M}{|t|} \right)^{\frac{1}{2}} \right)} \right)^{1-
\frac{2}{r}}
\end{array}
\label{Eqn:EstKMtf}
\end{equation}
and $r^{'}$ such that $\frac{1}{r}+\frac{1}{r^{'}}=1$. Observe that
if $(q,r)$ is wave admissible and $(q,r) \neq (\infty,2)$ then
$\frac{1}{q} + \frac{d}{2r} < \frac{d}{4} $. Therefore there are two
cases \vspace{2mm}

First we estimate $\| T_{1} T_{1}^{*} \|_{L_{t}^{q} L_{x}^{r}}$.
There are two cases

\begin{itemize}

\item \textbf{Case 1}: $r > 2$. Then since $(q,r)$ is wave admissible and $(q,r) \neq (\infty ,2)$ we also have
$\frac{1}{q} + \frac{d}{2r} < \frac{d}{4}$ and by (\ref{Eqn:kern1}),
Young's inequality and (\ref{Eqn:EstK1})

\begin{equation}
\begin{array}{ll}
\| T_{1} T_{1}^{\ast} g \|_{L_{t}^{q} L_{x}^{r}} & \lesssim \| \min{
\left(1, \frac{1}{|t|^{\frac{d}{2}}} \right)}^{1-\frac{2}{r}}
\|_{L_{t}^{\frac{q}{2}}} \| g \|_{L_{t}^{q^{'}} L_{x}^{r^{'}}} \\
& \lesssim \| g \|_{L_{t}^{q^{'}} L_{x}^{r^{'}}}
\end{array}
\end{equation}

\item \textbf{Case 2}: $r=2$. Then $q=\infty$. Then by (\ref{Eqn:kern1}) and (\ref{Eqn:PlK1}) we get (\ref{Eqn:ResT1}).

\end{itemize}

\vspace{2mm}

We turn to (\ref{Eqn:ResTM}). We write $\widetilde{K_{M}}=
\widetilde{K_{M,a}} + \widetilde{K_{M,b}} + \widetilde{K_{M,c}}$ in
(\ref{Eqn:IneqConv}) with $\widetilde{K_{M,a}}:= \widetilde{K_{M}}
\chi_{ _{ |t| \leq \frac{1}{M}}}$,
$\widetilde{K_{M,b}}:=\widetilde{K_{M}} \chi_{_{\frac{1}{M} \leq |t|
\leq M}}$ and $\widetilde{K_{M,c}}:= \widetilde{K_{M}} \chi_{_{|t|
\geq M}}$. We have by Young's inequality and (\ref{Eqn:qrscal})

\begin{equation}
\begin{array}{ll}
\left\| \widetilde{K_{M,a}}(t-t^{'}) \, \| g(t^{'},.)
\|_{L_{x}^{r^{'}}} \right\|_{L_{t}^{q}} & \lesssim M^{d \left(1
-\frac{2}{r} \right)}  \| \chi _{ |t| \leq \frac{1}{M}}
\|_{L_{t}^\frac{q}{2}} \| g \|_{L_{t}^{q^{'}}
L_{x}^{r^{'}}} \\
& \lesssim M^{2m} \| g \|_{L_{t}^{q^{'}} L_{x}^{r^{'}}}
\end{array}
\label{Eqn:EstKma}
\end{equation}
To estimate $\| \widetilde{K_{M,b}}(t-t^{'}) \ast \| g(t^{'},.)
\|_{L^{r^{'}}}  \|_{L_{t}^{q}}$ there are two cases

\begin{itemize}

\item \textbf{Case 1}: $\frac{1}{q} + \frac{d-1}{2r} < \frac{d-1}{4}$
By Young's inequality, (\ref{Eqn:EstKMtf}) and (\ref{Eqn:qrscal})
we have

\begin{equation}
\begin{array}{ll}
\left\| \widetilde{K_{M,b}}(t-t^{'}) \, \| g(t^{'},.)
\|_{L_{x}^{r^{'}}}   \right\|_{L_{t}^{q}} & \lesssim \|
\chi_{\frac{1}{M} \leq |t| \leq M} \frac{M^{d \left( 1- \frac{2}{r}
\right)}}{(M |t|)^{\frac{d-1}{2} \left( 1 - \frac{2}{r} \right)}}
\|_{L_{t}^{\frac{q}{2}}} \| g \|_{L_{t}^{q^{'}} L_{x}^{r^{'}}} \\
& \lesssim M^{2m} \| g \|_{L_{t}^{q^{'}}L_{x}^{r^{'}}}
\end{array}
\label{Eqn:EstKmb1}
\end{equation}

\item \textbf{Case 2}: $\frac{1}{q} + \frac{d-1}{2r} = \frac{d-1}{4}$. By (\ref{Eqn:EstKMtf}) we have

\begin{equation}
\begin{array}{ll}
\widetilde{K_{M,b}}(t-t^{'}) \, \| g(t^{'},.) \|_{L_{x}^{r^{'}}}  &
\lesssim M^{\frac{d+1}{2} \left( 1- \frac{2}{r} \right)}
 \frac{ \| g(t^{'},.) \|_{L^{r^{'}}}}{|t-t^{'}|^{\frac{d-1}{2} \left(1 -\frac{2}{r} \right)}}  \\
& \lesssim M^{d \left( 1 -\frac{2}{r} \right) -\frac{2}{q}}
 \frac{ \| g(t^{'},.)  \|_{L^{r^{'}}}}{|t-t^{'}|^{\frac{2}{q}}} \\
\end{array}
\end{equation}
By (\ref{Eqn:EdPt}), (\ref{Eqn:qrscal}) and Hardy-Littlewood-Sobolev
inequality \cite{stein}

\begin{equation}
\begin{array}{ll}
\left\| \widetilde{K_{M,b}}(t-t^{'}) \, \| g(t^{'},.)
\|_{L_{x}^{r^{'}}}  \right\|_{L_{t}^{q}}
 & \lesssim  M^{2m} \| g \|_{L_{t}^{q^{'}} L_{x}^{r^{'}}}
\end{array}
\label{Eqn:EstKmb2}
\end{equation}

\end{itemize}
We estimate $ \left\| \widetilde{K_{M,c}}(t-t^{'}) \, \| g(t^{'},.)
\|_{L^{r^{'}}}  \right\|_{L_{t}^{q}} $ by applying Young inequality,
(\ref{Eqn:qrscal}) and (\ref{Eqn:Waveadm}) i.e

\begin{equation}
\begin{array}{ll}
 \left\| \widetilde{K_{M,c}}(t-t^{'}) \, \| g(t^{'},.) \|_{L_{x}^{r^{'}}}  \right\|_{L_{t}^{q}} & \lesssim
 M^{ \left( \frac{d}{2} +1 \right) \left( 1- \frac{2}{r} \right)} \| \chi_{_{|t| \geq M}} \frac{1}{|t|^{\frac{d}{2}
 \left( 1- \frac{2}{r} \right) }} \|_{L_{t}^{\frac{q}{2}}} \| g \|_{L_{t}^{q'} L_{x}^{r'}} \\
& \lesssim M^{\frac{2}{q} + 1 -\frac{2}{r}} \| g \|_{L_{t}^{q^{'}} L_{x}^{r^{'}}} \\
& \lesssim M^{2m} \| g \|_{L_{t}^{q^{'}} L_{x}^{r^{'}}}
\end{array}
\label{Eqn:EstKmc}
\end{equation}
By (\ref{Eqn:IneqConv}), (\ref{Eqn:EstKma}), (\ref{Eqn:EstKmb1}),
(\ref{Eqn:EstKmb2}) and (\ref{Eqn:EstKmc}) we get (\ref{Eqn:ResTM}).

\subsection{Proof of (\ref{Eqn:HmNLin1})}

By decomposition and substitution it suffices to prove

\begin{equation}
\begin{array}{ll}
\| \int_{t^{'} < t} e^{i(t-t^{'}) \langle D \rangle} Q(t^{'}) \,
dt^{'} \|_{L_{t}^{q} L_{x}^{r}} & \lesssim \| \langle D \rangle  Q
\| _{L_{t}^{\widetilde{q}} L_{x}^{\widetilde{r}} }
\end{array}
\end{equation}
By Christ-Kisilev lemma \cite{SmSog} \footnote{an original proof of
this lemma can be found in \cite{chrkis}} it suffices in fact to
prove

\begin{equation}
\begin{array}{ll}
\| \int e^{i(t-t^{'}) \langle D \rangle} Q(t^{'}) \, dt^{'}
\|_{L_{t}^{q} L_{x}^{r}} & \lesssim \| \langle D \rangle Q \|
_{L_{t}^{\widetilde{q}} L_{x}^{\widetilde{r}} }
\end{array}
\label{Eqn:HmLin1pr}
\end{equation}
If we could prove

\begin{equation}
\begin{array}{ll}
\| \int e^{i(t-t^{'}) \langle D \rangle} P_{\leq 1} Q(t^{'}) \,
dt^{'} \|_{L_{t}^{q} L_{x}^{r}} & \lesssim \| Q \|
_{L_{t}^{\widetilde{q}} L_{x}^{\widetilde{r}}}
\end{array}
\label{Eqn:LwFreqInhom}
\end{equation}
and

\begin{equation}
\begin{array}{ll}
\| \int e^{i(t-t^{'}) \langle D \rangle } P_{M} Q(t^{'}) \, dt^{'}
\|_{L_{t}^{q} L_{x}^{r}} & \lesssim  M \| Q
\|_{L_{t}^{\widetilde{q}} L_{x}^{\widetilde{r}}}
\end{array}
\label{Eqn:HghFreqInhom}
\end{equation}
then (\ref{Eqn:HmLin1}) would follow. Indeed introducing
$\widetilde{P_{\leq 1}} $ and $ \widetilde{P_{M}} $  as in the
previous subsection we have

\begin{equation}
\begin{array}{ll}
\| \int e^{i(t-t^{'}) \langle D \rangle} P_{\leq 1} Q(t^{'}) \,
dt^{'} \|_{L_{t}^{q} L_{x}^{r}} & \lesssim \| \widetilde{P_{\leq 1}}
Q \|
_{L_{t}^{\widetilde{q}} L_{x}^{\widetilde{r}}} \\
& \lesssim \|  \langle D \rangle Q \|_{L_{t}^{\widetilde{q}}
L_{x}^{\widetilde{r}}}
\end{array}
\end{equation}
and

\begin{equation}
\begin{array}{ll}
\| \int e^{i(t-t^{'}) \langle D \rangle} P_{M} Q(t^{'}) \, dt^{'}
\|_{L_{t}^{q} L_{x}^{r}} & \lesssim M \| \widetilde{P_{M}} Q \|
_{L_{t}^{\widetilde{q}} L_{x}^{\widetilde{r}}} \\
& \lesssim \| \widetilde{P_{M}} \langle D \rangle Q
\|_{L_{t}^{\widetilde{q}} L_{x}^{\widetilde{r}}}
\end{array}
\end{equation}
Therefore we have 

\begin{equation}
\begin{array}{l}
\| \int e^{i(t-t^{'}) \langle D \rangle} P_{> 1} Q(t^{'}) \, dt^{'}
\|_{L_{t}^{q} L_{x}^{r}}  \lesssim \left\|  \left( \sum\limits_{  M
\in
2^{\mathbb{Z}}, \, M > 1}  | \int e^{i(t-t^{'}) \langle D \rangle} P_{M} Q(t{'}) \, dt^{'} |^{2} \right)^{\frac{1}{2}} \right\|_{L_{t}^{q} L_{x}^{r}}  \\
\lesssim \left( \sum \limits_{ M \in
2^{\mathbb{Z}}, \, M > 1}  \| \int e^{i(t-t^{'}) \langle D \rangle} P_{M} Q(t^{'}) \, dt^{'} \|^{2}_{L_{t}^{q} L_{x}^{r}} \right)^{\frac{1}{2}} \\
 \lesssim  \left( \sum \limits_{ M \in
2^{\mathbb{Z}}, \,  M > 1}  \| \widetilde{P_{M}} \langle D \rangle Q \|^{2}_{L_{t}^{\widetilde{q}} L_{x}^{\widetilde{r}}} \right)^{\frac{1}{2}} \\
 \lesssim  \left\|  \left( \sum \limits_{ M \in 2^{\mathbb{Z}}, \,
M > 1}  | \widetilde{P_{M}} \langle D \rangle Q|^{2} \right)^{\frac{1}{2}}
\right\|_{L_{t}^{\widetilde{q}}
L_{x}^{\widetilde{r}}} \\
 \lesssim \| \langle D \rangle Q \|_{L_{t}^{\widetilde{q}}
L_{x}^{\widetilde{r}}}
\end{array}
\label{Eqn:StepsInhom}
\end{equation}
Now we establish (\ref{Eqn:LwFreqInhom}). It is not difficult to see
from the proof of (\ref{Eqn:StrLocLin1}) and (\ref{Eqn:StrLocLinM})
that we also have

\begin{equation}
\begin{array}{ll}
\| e^{it \langle D \rangle} P^{\frac{1}{2}}_{\leq 4} f \|_{L_{t}^{q}
L_{x}^{r}} & \lesssim \| f \|_{L^{2}}
\end{array}
\label{Eqn:Lbarqr4}
\end{equation}
and

\begin{equation}
\begin{array}{ll}
\| e^{it \langle D \rangle} P^{\frac{1}{2}}_{\leq 1} f \|_{L_{t}^{q}
L_{x}^{r}} & \lesssim \| f \|_{L^{2}}
\end{array}
\label{Eqn:Lbarqr1}
\end{equation}
for every Schwartz function $f$. A dual statement of
(\ref{Eqn:Lbarqr1}) is

\begin{equation}
\begin{array}{ll}
\| \int e^{-i t^{'} \langle D \rangle} P^{\frac{1}{2}}_{ \leq 1}
Q(t^{'}) \, dt^{'} \|_{L^{2}} & \lesssim \| Q
\|_{L_{t}^{\widetilde{q}} L_{x}^{\widetilde{r}}}
\end{array}
\label{Eqn:LDualSt}
\end{equation}
Composing (\ref{Eqn:Lbarqr4}) with (\ref{Eqn:LDualSt}) we get (\ref{Eqn:LwFreqInhom}). \\
We turn to (\ref{Eqn:HghFreqInhom}). We need to prove

\begin{equation}
\begin{array}{ll}
\| \int e^{i(t-t^{'}) \langle \xi \rangle} \psi \left( \frac{\xi}{M}
\right) \widehat{Q}(t^{'},\xi) \, d t^{'} \, e^{i \xi \cdot x} \, d
\xi \|_{L_{t}^{q} L_{x}^{r}} & \lesssim M \| Q
\|_{L_{t}^{\widetilde{q}} L_{x}^{\widetilde{r}}}
\end{array}
\label{Eqn:ToProveExp}
\end{equation}
By the change of variable  $\left( \xi, \, t^{'}  \right)
\rightarrow \left( \frac{\xi}{M}, \, M t^{'} \right) $ we are
reduced showing

\begin{equation}
\begin{array}{ll}
\| \int e^{i(Mt-t^{'}) \left( |\xi|^{2} + \frac{1}{M^{2}}
\right)^{\frac{1}{2}}} \psi(\xi) \widehat{Q \left( \frac{t^{'}}{M},
\, \frac{.}{M} \right)}(\xi)  \, dt^{'} \, e^{i M x \cdot \xi} \, d
\xi \|_{L_{t}^{q} L_{x}^{r}} & \lesssim M^{2} \| Q
\|_{L_{t}^{\widetilde{q}} L_{x}^{\widetilde{r}}}
\end{array}
\label{Eqn:IneqMInterm}
\end{equation}
If we could prove that for every Schwartz function $G$

\begin{equation}
\begin{array}{ll}
\| S_{M} G \|_{L_{t}^{q} L_{x}^{r}} & \lesssim \| G
\|_{L_{t}^{\widetilde{q}} L_{x}^{\widetilde{r}}}
\end{array}
\label{Eqn:SMIneq}
\end{equation}
with
\begin{equation}
\begin{array}{ll}
S_{M} G & := \int e^{i(t-t^{'}) \left( |\xi|^{2} + \frac{1}{M^{2}}
\right)^{\frac{1}{2}}} \psi(\xi) \widehat{G}(t^{'},\xi) \, dt^{'}
\, e^{i \xi \cdot x} d\xi
\end{array}
\label{Eqn:DefSM}
\end{equation}
then (\ref{Eqn:IneqMInterm}) would hold. Indeed by
(\ref{Eqn:ScalingInhomqtildeq}) we have

\begin{equation}
\begin{array}{ll}
\| \int e^{i(Mt-t^{'}) \left( |\xi|^{2} + \frac{1}{M^{2}}
\right)^{\frac{1}{2}}} \psi(\xi) \widehat{Q \left( \frac{t^{'}}{M},
\, \frac{.}{M} \right)}(\xi)  \, dt^{'} \, e^{i M x \cdot \xi} \, d
\xi \|_{L_{t}^{q} L_{x}^{r}} \\  = \| S_{M} \left( Q \left(
\frac{.}{M},\frac{.}{M} \right) \right)(Mt, \, Mx) \|_{L_{t}^{q} L_{x}^{r}} \\
 \lesssim M^{\frac{1}{\widetilde{q}} + \frac{d}{\widetilde{r}} - \frac{1}{q} -\frac{d}{r}} \| Q \|_{L_{t}^{\widetilde{q}}
L_{x}^{\widetilde{r}}} \\
 \lesssim M^{2} \| Q \|_{L_{t}^{\widetilde{q}} L_{x}^{\widetilde{r}}}
\end{array}
\label{Eqn:ScalingSM}
\end{equation}
By duality and composition with (\ref{Eqn:Lbarqr4}) it suffices to
show

\begin{equation}
\begin{array}{ll}
\| e^{it \left( D^{2} + \frac{1}{M^{2}} \right)^{\frac{1}{2}}}
P^{\frac{1}{2}}_{1} f  \|_{L_{t}^{\tilde{q}^{'}} L_{x}^{\tilde{r}^{'}}} & \lesssim \| f
\|_{L^{2}}
\end{array}
\label{Eqn:RemainInt}
\end{equation}
Again it is not difficult to see from the proof of (\ref{Eqn:TMEst})
that

\begin{equation}
\begin{array}{ll}
\| e^{it \langle D \rangle} P^{\frac{1}{2}}_{M} f \|_{L_{t}^{\tilde{q}^{'}}
L_{x}^{\tilde{r}^{'}}} & \lesssim \| f \|_{H^{1-m}}
\end{array}
\label{Eqn:TMEstvar1}
\end{equation}
But after performing the change of variable $\xi \rightarrow M \xi$
we have by (\ref{Eqn:qrscal}) and (\ref{Eqn:TMEstvar1})

\begin{equation}
\begin{array}{ll}
\| e^{it \left( D^{2} + \frac{1}{M^{2}} \right)^{\frac{1}{2}}}
P^{\frac{1}{2}}_{1} f  \|_{L_{t}^{\tilde{q}^{'}} L_{x}^{\tilde{r}^{'}}} & =  \| \int e^{i
\frac{t}{M} \left( |\xi|^{2} + 1 \right)^{\frac{1}{2}}}
\psi^{\frac{1}{2}} \left( \frac{\xi}{M}   \right) \widehat{f(M \,
.)}(\xi) e^{i \frac{x}{M} \cdot \xi} \, d \xi
\|_{L_{t}^{\tilde{q}^{'}} L_{x}^{\tilde{r}^{'}}} \\
& = \| \left(  e^{it \langle D \rangle} P^{\frac{1}{2}}_{M} \right) \left( \widetilde{P_{M}}^{\frac{1}{2}} f (M.) \right) \left( \frac{t}{M}, \frac{x}{M} \right) \|_{L_{t}^{\tilde{q}^{'}} L_{x}^{\tilde{r}^{'}}} \\
& \lesssim M^{\frac{1}{\tilde{q}^{'}} + \frac{d}{\tilde{r}^{'}}} \| \widetilde{P_{M}}^{\frac{1}{2}}  f (M \, .)  \|_{H^{1-m}} \\
& \lesssim M^{\frac{1}{\tilde{q}^{'}} + \frac{d}{\tilde{r}^{'}} -\frac{d}{2} + 1-m}\| f \|_{L^{2}} \\
& \lesssim \| f \|_{L^{2}}
\end{array}
\label{Eqn:ChVar}
\end{equation}

\subsection{ Proof of (\ref{Eqn:HmNLin2})}

By decomposition, substitution and Christ-Kisilev lemma \cite{SmSog}
it suffices to prove

\begin{equation}
\begin{array}{ll}
\| \int  e^{i (t- t^{'}) \langle D \rangle} Q  \|_{L_{t}^{\infty}
L_{x}^{2}} & \lesssim  \| \langle D \rangle^{1-m} Q
\|_{L_{t}^{\widetilde{q}} L_{x}^{\widetilde{r}}}
\end{array}
\label{Eqn:HNLin2ToPr}
\end{equation}
If we could prove

\begin{equation}
\begin{array}{ll}
\| \int  e^{i (t- t^{'}) \langle D \rangle} P_{\leq 1} Q(t^{'}) \, d
t^{'} \|_{L_{t}^{\infty} L_{x}^{2}} & \lesssim \| Q
\|_{L_{t}^{\tilde{q}} L_{x}^{\tilde{r}}}
\end{array}
\label{Eqn:1Toprove}
\end{equation}
and

\begin{equation}
\begin{array}{ll}
\| \int e^{i(t - t^{'}) \langle D \rangle} P_{M} Q(t^{'}) \, dt^{'}
\|_{L_{t}^{\infty} L_{x}^{2}} & \lesssim M^{1-m} \| Q
\|_{L_{t}^{\tilde{q}} L_{x}^{\tilde{r}}}
\end{array}
\label{Eqn:MToprove}
\end{equation}
then (\ref{Eqn:HNLin2ToPr}) would follow. Indeed

\begin{equation}
\begin{array}{ll}
\| \int e^{i(t-t^{'}) \langle D \rangle} P_{\leq 1} Q(t^{'}) \,
dt^{'} \|_{L_{t}^{\infty} L_{x}^{2}} & \lesssim \|
\widetilde{P_{\leq 1}} Q
\|_{L_{t}^{\tilde{q}} L_{x}^{\tilde{r}}} \\
& \lesssim  \| \langle D \rangle^{1-m} Q \|_{L_{t}^{\tilde{q}}
L_{x}^{\tilde{r}}}
\end{array}
\end{equation}
and

\begin{equation}
\begin{array}{ll}
\| \int e^{i(t-t^{'}) \langle D \rangle} P_{M} Q(t^{'}) \, dt^{'}
\|_{L_{t}^{\infty} L_{x}^{2}}
& \lesssim M^{1-m} \| \widetilde{P_{M}} Q  \|_{L_{t}^{\tilde{q}} L_{x}^{\tilde{r}}} \\
& \lesssim \| \widetilde{P_{M}} \langle D \rangle^{1-m} Q
\|_{L_{t}^{\tilde{q}} L_{x}^{\tilde{r}}}
\end{array}
\end{equation}
Therefore following the same steps to those in
(\ref{Eqn:StepsInhom}) we get (\ref{Eqn:HNLin2ToPr}).

(\ref{Eqn:1Toprove}) follows from the composition of the trivial
inequality  $ \| e^{it \langle D \rangle} P^{\frac{1}{2}}_{\leq 1} f
\|_{L_{t}^{\infty} L_{x}^{2}} \lesssim  \| f \|_{L^{2}} $ and
(\ref{Eqn:LDualSt}).

We turn to (\ref{Eqn:MToprove}). We need to prove

\begin{equation}
\begin{array}{ll}
\| \int e^{i(t-t^{'}) \langle \xi \rangle} \psi \left( \frac{\xi}{M}
\right) \widehat{Q}(t^{'},\xi) \, d t^{'} \, e^{i \xi \cdot x} \, d
\xi \|_{L_{t}^{\infty} L_{x}^{2}} & \lesssim M^{1-m} \| Q
\|_{L_{t}^{\tilde{q}} L_{x}^{\tilde{r}}}
\end{array}
\end{equation}
Again by the change of variable  $\left( \xi, \, t^{'}  \right)
\rightarrow \left( \frac{\xi}{M}, \, M t^{'} \right) $ it suffices
to show

\begin{equation}
\begin{array}{ll}
\| \int e^{i(Mt-t^{'}) \left( |\xi|^{2} + \frac{1}{M^{2}}
\right)^{\frac{1}{2}}} \psi(\xi) \widehat{Q \left( \frac{t^{'}}{M},
\, \frac{.}{M} \right)}(\xi)  \, dt^{'} \, e^{i M x \cdot \xi} \, d
\xi \|_{L_{t}^{\infty} L_{x}^{2}} \\
\lesssim M^{2-m} \| Q
\|_{L_{t}^{\tilde{q}} L_{x}^{\tilde{r}}}
\end{array}
\label{Eqn:IneqMIntermLast}
\end{equation}
If we could prove for any Schwartz function

\begin{equation}
\begin{array}{ll}
\| S_{M} G \|_{L_{t}^{\infty} L_{x}^{2}} & \lesssim \| G
\|_{L_{t}^{\tilde{q}} L_{x}^{\tilde{r}}}
\end{array}
\label{Eqn:ToproveSminfty}
\end{equation}
with $S_{M}$ defined in (\ref{Eqn:DefSM}) then substituting $q$, $r$
for $\infty$, $2$ respectively in (\ref{Eqn:ScalingSM}) we have

\begin{equation}
\begin{array}{l}
\| \int e^{i(Mt-t^{'}) \left( |\xi|^{2} + \frac{1}{M^{2}}
\right)^{\frac{1}{2}}} \psi(\xi) \widehat{Q \left( \frac{t^{'}}{M},
\, \frac{.}{M} \right)}(\xi)  \, dt^{'} \, e^{i M x \cdot \xi} \, d
\xi \|_{L_{t}^{\infty} L_{x}^{2}} \\
\lesssim M^{\frac{1}{\tilde{q}}
+ \frac{d}{\tilde{r}} -\frac{d}{2}} \| Q \|_{L_{t}^{\tilde{q}}
L_{x}^{\tilde{r}}} \\
\lesssim M^{2-m} \| Q \|_{L_{t}^{\tilde{q}} L_{x}^{\tilde{r}}}
\end{array}
\end{equation}
where in the last inequality we used (\ref{Eqn:qrscal}) and
(\ref{Eqn:ScalingInhomqtildeq}). It remains to prove
(\ref{Eqn:ToproveSminfty}). By duality and composition with the
trivial inequality $\| e^{i t \langle D \rangle} P^{\frac{1}{2}}_{\leq 4} f
\|_{L_{t}^{\infty} L_{x}^{2}} \lesssim  \| f \|_{L^{2}}$ it suffices
to show (\ref{Eqn:RemainInt}), which has already been established.

\end{document}